\documentclass[12pt,leqno,fleqn,epsfig]{article}
\usepackage{amssymb, epsfig, amsmath, amsthm}
\usepackage{mathrsfs} 

\newcommand{\R}{\mathbb{R}}
\newcommand{\PP}{\mathbb{P}}

\newcommand{\N}{\mathbb{N}}

\newcommand{\curl}{\operatorname{curl}}

\newcommand{\bb}{\begin{equation}}
\newcommand{\ee}{\end{equation}}
\newcommand{\bq}{\begin{eqnarray}}
\newcommand{\eq}{\end{eqnarray}}
\newcommand{\bqn}{\begin{eqnarray*}}
\newcommand{\eqn}{\end{eqnarray*}}

\newcommand{\var}{\varepsilon}
\renewcommand{\o}{\omega}

\newcommand{\intl}{\int\limits}

\newcommand{\Beweisende}{\rule{0.2cm}{0.2cm}}

\newcommand{\D}{\displaystyle}

\newcommand{\intmw}{{\int\hspace{-830000sp}-\!\!}}

\newcounter{secnum}

\newtheorem{thm}{Theorem}[section]
\newtheorem{cor}[thm]{Corollary}
\newtheorem{lem}[thm]{Lemma}

\theoremstyle{definition}

\newtheorem{rem}[thm]{Remark}

\title{Removing  Type II singularities off the axis for  the 3D axisymmetric Euler equations} 
 
\author{Dongho Chae$^*$  and J\"{o}rg Wolf $^\dagger$\\
\ \\
Department of Mathematics\\
Chung-Ang University\\
 Seoul 156-756, Republic of Korea\\
  ($*$)e-mail: dchae@cau.ac.kr\\
($\dagger$)e-mail: jwolf2603@cau.ac.kr}
\date{}
\begin{document}
\maketitle
\begin{abstract}
In this paper we obtain new  local blow-up criterion for smooth  axisymmetric solutions to the 3D incompressible Euler equation.  If the vorticity satisfies $ \intl_{0}^{t_*} (t_*-t)   \| \omega   (t)\|_{ L^\infty(B(x_{ \ast}, R_0))}   dt <+\infty$ for a ball $B(x_{ \ast}, R_0)$ away from the axis of symmetry, then there exists no singularity  at $t=t_*$ in the torus $T(x_*, R)$ 
generated by rotation of the ball $B(x_{ \ast}, R_0)$ around the axis. This implies that  possible singularity  at $t=t_*$ in  the torus  $T(x_*, R)$  is excluded if the vorticity satisfies the  blow-up rate $ \|\o (t)\|_{L^\infty (T(x_*, R))}= O\left(\frac{1}{(t_*-t)^\gamma}\right)$  as $t\to t_*$, where  $\gamma <2$,  and the torus $T(x_*, R)$ does not touch the axis.
  \noindent{\bf
keywords:} Euler equations, axisymmetric solution,  local blow-up criterion, Type II blow-up

\end{abstract}

\section{Introduction}
\label{sec:-1}
\setcounter{secnum}{\value{section} \setcounter{equation}{0}
\renewcommand{\theequation}{\mbox{\arabic{secnum}.\arabic{equation}}}}

We consider the  Euler equations in the domain $ \R^{3}\times (0, +\infty)$
\begin{equation}
\begin{cases}
\partial _t v  + (v \cdot \nabla ) v  = -\nabla p,\\
\quad\nabla \cdot  v =0,\\
v(x,0)=v_0(x)\quad  \forall  x\in\R^{3}
\end{cases}
\label{euler}
\end{equation}
where $v=(v_1 (x,t), v_2 (x,t), v_3(x,t))$ is the velocity of the fluid, and $p=p(x,t)$ represents the pressure. Here, $ v_0$ stands for  initial velocity.  Given  $v_0 \in W^{k,q} (\Bbb R^3)$, $k>3/q+1$, $1<q<+\infty$,  the local in time well-posedness is well-known due to Kato and Ponce\cite{kat}. The question of spontaneous apparition of singularity in  finite time, however, is an outstanding open problem in the mathematical fluid mechanics(see e.g.\cite{maj, con, bar} for surveys  of  the problem and the related results).  There are also many numerical experiments on the problem(see e.g. \cite{pum, ker, gra, ohk, bre, hou}).

\vspace{0.2cm}
\hspace{0.5cm}
We say a local in time smooth solution $v\in C([0, t_*); W^{k,q}(\Bbb R^3))$, $k>3/q+1$, $1<q<+\infty$,
 blows up (or equivalently becomes singular) at $t=t_*$ if 
\begin{equation}\label{noblow}
\limsup_{t\nearrow t_*} \|v(t)\|_{W^{k,q}(\Bbb R^3)}=+\infty \qquad \forall k\ge \frac{3}{q}+1.
\end{equation}
The  celebrated Beale-Kato-Majda(BKM) criterion\cite{bea} says that \eqref{noblow} holds if and only if 
\begin{equation}\label{bk}
\int_{0} ^{t_*} \|\omega (t)\|_{L^\infty(\Bbb R^3)} dt =+\infty, \qquad \omega =\nabla \times v.
\end{equation} 
See also \cite{con, den} for geometric type criteria.  Later, Kozono and Taniuchi\cite{koz} improved \eqref{bk}, replacing $\|\omega(t)\|_{L^\infty(\Bbb R^3)}$ in \eqref{bk} by a weaker norm $\|\omega(t)\|_{BMO(\Bbb R^3)}$. In a very recent paper authors of current paper  obtained a localized version of the criterions of \cite{bea, koz}, which says that in order to check the blow-up  for the  solution at particular space-time point $(x_*,t_*)\in \Bbb R^3\times (0, +\infty)$, it suffices to check if there exists $r>0$ such that  
\bb\int_0^{t_*} \|\o(t)\|_{L^\infty(B(x_*,r))}dt=+\infty,
\label{bk1}
\ee
where and hereafter we use the notation  $B(x, r)=\{ y\in \Bbb R^n\, |\, |x-y|<r\}$, and $B(r):=B(0, r)$ (In fact, we derived the localized criterion with $\|\o(t)\|_{L^\infty(B(x_*,r))}$ of  \eqref{bk1}
replaced by a weaker norm $\|\omega (t)\|_{BMO(B(x_*,r))}$).  Note that the local BKM criterion does not rule out the scenario of the blow-up rate 
\bb
\label{rate1}
\lim\sup_{t\to t_*}(t_*-t)^\gamma \|\o(t)\|_{L^\infty(B(x_*,r))}<+\infty
\ee
for some $x_*\in \Bbb R^3, r>0$, and for  $\gamma \ge 1$,  which includes the case
\bb
\label{rate2}
\|\o(t)\|_{L^\infty(B(x_*,r))} =O \left(\frac{1}{(t_* -t)^\gamma }\right) \quad \text{as}\quad t\to t_* .
\ee

\vspace{0.2cm}
\hspace{0.5cm}
In this paper we study the finite time blow-up problem of system \eqref{euler} under the assumption of axial symmetry. In this case also there are many previous studies from theoretical(e.g. \cite{cha4, cha2, cha3a, cha3b, cafl, dan})  and  computational(see \cite{hou, pum} and the references therein) aspects respectively.
Our  purpose  is to obtain a  local BKM criterion \eqref{bk} for the axisymmetric Euler equation off the axis region, which implies that  we can rule out the singularity having the blow-up rate including 
\eqref{rate1} and \eqref{rate2} with $\gamma<2$ at the region.  We say the vorticity $ \o$ 
has  Type I blow-up at $ t_*$  in a domain $\mathcal{D}\subset \Bbb R^3$ if 
\begin{equation}
\lim\sup_{t\to t_*} (t_*-t)\| \o(t)\|_{ \infty(\mathcal{D})} <+\infty. 
\label{type1}
\end{equation}
 Otherwise, we say it is of Type II.   For comparison we recall the status of the Type I blow-up problem in the axisymmetric Navier-Stokes equations. Due to the the well-known result of Caffarelli-Kohn-Nirenberg \cite{caf} for a suitable weak solutions  to the 3D Navier-Stokes equations  the one-dimensional Hausdorff measure of the possible singular set must be zero, and therefore all possible  singularities should be on the axis of symmetry. On the other hand, according  to the works of \cite{ser} and \cite{che1} independently there exists no Type I singularity
 on the axis of symmetry. For the case of our axisymmetric Euler equations there exists no available partial regularity type results similar to the Navier-Stokes equations,  and therefore we cannot rule out singularity off the axis by an simple argument. 
 Our main result, Theorem \ref{thm1.1} (and its immediate consequence Theorem \ref{thm1.2}) shows that  there exists no Type I 
 (and also part of Type II) singularity off the axis region. 
 
\vspace{0.2cm}
\hspace{0.5cm}
Below we briefly introduce the notion of axisymmetric flow,  which is found in standard literature(e.g. \cite{maj}).  We 
say $v$ is an axisymmetric solution of the Euler equations if  $v$ solves  \eqref{euler}, and can be written as
\[
v= v^ r (r,x_3,t) e_r + v^{ \theta }(r, x_3,t)e_\theta   +v^3 (r,x_3,t) e_3, 
\]   
 where 
 \[
e_r = (\frac{x_1}{r} ,\frac{x_2}{r}, 0),\quad  e_\theta  = ( \frac{x_2}{r},\frac{-x_1}{r}, 0), \quad  e_3 =  (0, 0, 1),\quad  r = \sqrt{x_1^2 + x_2^2}.
\]
are the canonical basis of the cylindrical coordinate system.
The  Euler equations \eqref{euler} for an axisymmetric solution  turn into the following equations.
\begin{align}
\partial _t v^r + v^r \partial _r v^r +v^3\partial _3 v^r &= - \partial _r p+ \frac{(v^\theta) ^2}{r}, 
\label{eu-r}
\\
\partial _t v^\theta  + v^r \partial _r v^\theta  +v^3  \partial _3 v^\theta  &= - \frac{v^r v^\theta }{r}, 
\label{eu-th}
\\
\partial _t v^3 + v^r  \partial _r v^3 +v^3 \partial _3 v^3 &= - \partial _3 p,
\label{eu-3}
\\
\partial _r (rv^r) +\partial _3(rv^3)=0. 
\label{div}\end{align} 
Multiplying \eqref{eu-th} by $ r$,  we see that $ r v^\theta $ satisfies  the transport equation
\begin{equation}
\partial _t (r v^\theta) + v^r \partial _r (r v^\theta)+ v^3 \partial _3 (r v^\theta) =0. 
\label{trans}
\end{equation} 
For the vorticity $ \omega=\nabla \times v$ we get
\[
\omega = \omega ^r e_r +  \omega ^\theta  e_\theta  + \omega ^3 e_3, 
\]
where
\[
\omega ^r = - \partial _3 v^\theta,\quad  \omega ^\theta = \partial _3 v ^r - \partial_r v^3,\quad  \omega ^3= \frac{v^\theta }{r}+ \partial _r v^\theta.
\]
Applying $ \partial _3 $ to \eqref{eu-r} and applying $ \partial _r$ to \eqref{eu-3}, and taking the difference of the two equations, we obtain 
  the following equation for $ \omega ^\theta $
  \begin{equation}
  \partial _t \omega ^\theta + v^r  \partial _r \omega ^\theta  +v^3  \partial _3 \omega ^\theta=
  \frac{v^r \omega ^\theta }{r} + \partial _3 \frac{(v^\theta )^2}{r}. 
  \label{vort-th}
  \end{equation}
This leads to the equation 
  \begin{equation}
  \partial _t  \big(\frac{\omega ^\theta}{r} \big)+ v^r \partial _r \big(\frac{\omega ^\theta}{r} \big) +v^3  \partial _3 \big(\frac{\omega ^\theta}{r} \big)=
\frac{ \partial _3 (v^\theta )^2}{r^2}. 
  \label{vort-th1}
  \end{equation}

  Our first main theorem is the following improvement of the BKM theorem off the axis region.
\begin{thm}
\label{thm1.1}
Let $ v\in C([0,t_*); W^{2,\, q}(\R^{3})) \cap  L^\infty(0, t_*; L^{2}(\R^{3})),$ 
$ 3<q<+\infty, $  be an axisymmetric solution to  \eqref{euler} in $ \R^{3}\times (0,t_*)$. If the following condition is  fulfilled 
\begin{equation}
 \intl_{0}^{t_*} (t_*-t)   \| \omega   (t)\|_{ L^\infty(B(x_{ \ast}, R_0))}   dt <+\infty, \quad 
\label{type1i}
\end{equation}
for some ball    $B(x_{ \ast}, R_0) \subset \{ x \in \R^{3}\, |\, x_1^2+x_2^2 >0\}$, where $ \omega = \nabla \times v$. 
Then for all $ 0<R<R_0$ it holds  $ v\in C([0, t_*], W^{2,\, q}(B(x_{ \ast}, R)))$. In particular, this implies 
$v\in C([0, t_*], W^{2,\, q}(T(x_{ \ast}, R)))$. Here,  $ T(x_{ \ast}, R)$ stands for the torus generated by rotation of $ B(x_{ \ast}, R_0)$ around the axis, i.e. 
\[
T(x_{ \ast}, R) = \{x\in \R^{3} \,|\, \Big(\,\sqrt{ x_1^2+ x_2 ^2} - \rho_{ \ast} \Big)^2 +(x_3-x_{ 3, \ast})^2 < R^2 \},    
\]
where $ \rho_{ \ast} = \,\sqrt{x_{ 1, \ast}^2+ x_{2, \ast}^2}$. 

\end{thm}  

As an immediate consequence  of  this theorem we remove  some of Type II as well as Type I  singularities in terms of the vorticity blow-up rate off the axis. We have the following:

\begin{thm}
\label{thm1.2}
Let $ v\in C([0,t_*); W^{2,\, q}(\R^{3})) \cap  L^\infty(0,t_*; L^{2}(\R^{3})),$ 
$ 3<q<+\infty, $  be an axisymmetric solution to  \eqref{euler} in $ \R^{3}\times (0,t_*)$. 
Suppose the following vorticity blow-up rate condition holds
\begin{equation}
\sup_{ t\in (0,t_*)}(t_*-t)^2  \left[\log\left(\frac{1}{t_*-t}\right)\right]^{ \alpha }\| \omega (t)\|_{L^ \infty(B(x_{ \ast}, R_0))} < +\infty
\label{type2}
\end{equation}
for some   $ \alpha >1$ and some ball $B(x_{ \ast}, R_0) \subset \{ x \in \R^{3}\, |\, x_1^2+x_2^2 >0\}$. 
Then  $ v\in C([0, t_*]; W^{2,\, q}(T(x_{ \ast}, R))$ for all $ 0<R<R_0$. 
\end{thm}  
\begin{rem} 
 In particular, Theorem \ref{thm1.2} says that there exists no singularity at $t=t_*$  in the off-the axis region if  the the vorticity blow-up rate satisfies 
\bb\label{rem1}
\|\o (t)\|_{L^\infty (\Bbb R^3)}= O\left(\frac{1}{(t_*-t)^\gamma}\right), 
\ee
 as $t\to t_*$ if $1\leq \gamma <2$.  Due to the global BKM criterion, however, the singularity  in this case should happen only on the axis. It would be  interesting to compare this result with Tao's construction of a singular solution(see \cite[Fig. 3, pp.18]{tao}) for a {\em modified 
Euler system}, where $\gamma=1$  and the set of singularity  is  a circle around the axis.  
\end{rem}

\begin{rem}  In a recent numerical study of the blow-up  of the  axisymmetric  Euler equations by Luo and Hou\cite[pp.1766]{hou} they computed $\gamma= 2.45$ in \eqref{rem1}, which is consistent with our rigorous result, since their blow-up region is away from the axis, and near the boundary of the cylinder. As far as the authors know,  this is the first explicit computational result with $\gamma \geq 2$ for singularity in the axisymmetric Euler equations.
\end{rem}
The key observation in the proof of Theorem \ref{thm1.1} is that off the axis region, as is  well-known, the 3D Euler system behaves like a solution to the  2D Boussinesq system, the scaling property of which  is different from the 3D Euler equations, and our careful local analysis takes 
 full advantage of this property.  For the proof  we introduce  a new iteration scheme, the corresponding iteration lemma of which  is proved in the appedix B, and also we use a local version of the logarithmic Sobolev inequality proved  in \cite{cha5}.
 
 \vspace{0.2cm}
\hspace{0.5cm} 
    We note also that our local regularity criterion off the axis region,  is not an immediate consequence of the corresponding local criterion for the 2D Boussinesq system.  We need to  show further that $\int_0 ^{t_*} \|v(t)\|_{L^\infty (B(r, x_*))} dt <+\infty$ is satisfied under our hypothesis for the 3D axisymmetric case, which is accomplished by using 
 the estimate of extra component of vorticity, $\omega^\theta$ together with the local version of the Calder\'on-Zygmund inequality.

\vspace{0.2cm}
\hspace{0.5cm} 
 At this moment we do not know whether  Theorem\,\ref{thm1.1} continues to hold if we impose the condition 
\eqref{type1i} only for two components of vorticity 
$$\widetilde{\omega } := -\partial _3 v^\theta e_r + (\partial_r v^\theta + \frac{v^\theta}{r} )e_3.
$$  However,
a similar statement of Theorem\,\ref{thm1.1} holds,  if we replace the condition on $ \omega ^\theta $ by a corresponding condition on 
$ v^r$. More precisely, we have the following:  

\begin{thm}
\label{thm1.3}
Let $ v\in C([0,t_*); W^{2,\, q}(\R^{3})) \cap  L^\infty(0,t_*; L^{2}(\R^{3})),$ 
$ 2<q<+\infty, $  be an axisymmetric solution to  \eqref{euler} in $ \R^{3}\times (0,t_*)$. If the following two conditions are  fulfilled 
\begin{equation}
 \intl_{0}^{t_*} (t_*-t)   \| \widetilde{\omega }     (t)\|_{ L^\infty(\{ r>R\})}   dt <+\infty,  \quad  
 \intl_{0}^{t_*}  \| v^r    (t)\|_{ L^\infty(\{ r>R\})}   dt < +\infty
\label{type2i}
\end{equation}
for all   $ 0< R< +\infty$.  Then  $ v\in C([0, t_*]; W^{2,\, q}_{ loc}(\R^{3}  \setminus \{ r=0\}))$. 
\end{thm}  

\begin{rem}
Observing that in view of \eqref{trans}  $ r v^\theta $ is bounded in $ \R^{2}_+$ due to its conservation along the particle  trajectories generated by 
$$\tilde{v}= v^re_r+v^3 e_3,
$$
 it is immediately clear that \eqref{type2i}$ _1$  can be replaced  by
\begin{equation}
 \intl_{0}^{t_*} (t_*-t)   \| \nabla v^\theta    (t)\|_{ L^\infty(\{ r>R\})}   dt <+\infty. 
\label{type2j}
\end{equation} 
\end{rem}
The organization of this paper is as follows:

\vspace{0.2cm}
\hspace{0.5cm}
 In Section 2 we introduce a generalized 2D Boussinesq system, which includes the 3D axisymmetric Euler system off the axis and the standard 2D Boussinesq system as special cases. For such system we prove a local blow-up criterion,  Theorem \ref{thm3.1}, where certain integrability condition of `temperature' gradient  together with  the velocity integrability 
leads to a local non blow-up.  Establishing this theorem  is a major part for the proof the above main theorems. To prove Theorem \ref{thm3.1} this we  transform the equations from the generalized Boussineq system for $(u, \theta, w)$ into a new system for $(U, \Theta, \Omega)$. In order to deduce $W^{2,q}_{loc}$, $q>3$, estimates for the transformed functions $(U, \Theta)$, we need to handle differential inequality for integrals on different balls in the left and the right hand sides. Appearance of these different balls is originated from the use of cut-off functions, necessary for localizations, which is similar to the case of  deriving the Caccioppolli type inequalities in the elliptic equations.  In this case we cannot use  classical Gronwall's lemma to close the differential inequality. To overcome this difficulty we introduce a new of type iteration scheme, and  close the differential inequality by iterating a suitable sequence of differential inequalities. 

\vspace{0.2cm}
\hspace{0.5cm}
 In Section 3, using Theorem \ref{thm3.1}  of Section 2, we prove Theorem \ref{thm1.1}, where we only need to verify the integrability condition of the velocity of Theorem \ref{thm3.1} is satisfied. Using the fact that  the vorticity of the 3D axisymmetric  flow has an extra component  $\omega ^\theta$ other than $\tilde{\o}$,   the integrability of the velocity part can be proved. The proof of Theorem \ref{thm1.2} is an easy consequence of Theorem \ref{thm1.1}.  
 
 \vspace{0.2cm}
\hspace{0.5cm}
 In Section 4 
we prove Theorem \ref{thm1.3}, where the integrability condition of the radial velocity \eqref{type2i}$_2$ allows us to introduce the functional transform similarly to Section 2.  Using this transform, 
one can show that the integrability condition of $\nabla v^\theta$ implies the desired estimate of $\o^{\theta}$.  

\vspace{0.2cm}
\hspace{0.5cm}
In Appendix A we prove various  Gagliardo-Nirenberg type interpolation inequalities involving cut-off functions, which are necessary for the proof in Section 2. These include an improvement to local $BMO$ norm from  $L^\infty$ norm in  the known inequalities. 

\vspace{0.2cm}
\hspace{0.5cm}
In Appendix B we prove a Gronwall type  iteration lemma, which is crucial   to complete our  iteration scheme in  the proof of Theorem \ref{thm3.1} in Section 2.\\

{\em For simplicity of presentation and its proof below we consider the Euler system on the time interval $(-1,0)$, where $t=-1$ is the initial time and $t_*=0$ is the possible blow-up time.}

\section{Local BKM type criterion  for a generalized 2D Boussinesq system}
\label{sec:-3}
\setcounter{secnum}{\value{section} \setcounter{equation}{0}
\renewcommand{\theequation}{\mbox{\arabic{secnum}.\arabic{equation}}}}

The aim of this section is to prove a local regularity criterion of a generalized 2D Boussinesq equations, which 
include the 3D axisymmetric system as a special case.
Let  $B(1)$ be the unit ball in 
$\Bbb R^2$. 
We consider the following system on $ B(1)\times (-1, 0)$,
\begin{align}
\partial _t \theta  + a(y)u\cdot \nabla \theta &=0,
\label{bou1}
\\
  \partial _t  w   + a(y)u\cdot \nabla  w   &= b(y) \cdot \nabla  \theta, 
\label{bou2}
\\
\nabla \cdot u &=0,
\label{bou3}
\end{align}      
where $ u= (u_1(y_1, y_2, t), u_2(y_1, y_2, t)), \theta = \theta (y_1, y_2, t)$, and  $ u$ and $  w  $ are related by 
\begin{align}
  \curl u  &:= \partial_2 u_1 -\partial_1 u_2=d(y)  w  + e(y) \cdot u. 
\label{bou4}
\end{align}            
Here, $ a,b, d$ and  $e$ are coefficients in $C^2(\overline{B(1)} )$.  Note that the system  \eqref{trans}, \eqref{vort-th1} and \eqref{div}
reduces to the system   \eqref{bou1}-\eqref{bou4} with  identification   $ (r, x_3):=(y_1+2, y_2)$, $ u:= (rv^r, rv^3)$,  $ \theta := (rv^\theta )^2$,   $ w:= \frac{\omega ^\theta }{r}$ for the coefficient functions
  $a(y)=1/(y_1+2)$, $b(y)=(0,1/(y_1+2)^4)$ ,  $d(y)= (y_1+2)^2,$ and $ e(x)=(0, -1/(y_1+2))$.

Clearly, the above system covers the vorticity formulation of the standard  Boussinesq system in $ \R^{2}$. Namely, if $ (u, p, \theta )$ solves
the Boussinesq system in $ \R^{2}$, then $ (u,  w  , \theta )$ with $w=\partial_2 v_1-\partial _1 v_2$  solves \eqref{bou1}-\eqref{bou4} with 
$ a =d\equiv 1, b \equiv - e_2, $ and $ e \equiv 0$. Thus, the following result also applies to  the 2D Boussinesq system. 
Sufficient   conditions  for global regularity of solutions  to the  Boussinesq  system in whole $ \R^{2}$ has been proved  in \cite{cha6}.  Here we present the following local regularity condition  for the  system  \eqref{bou1}-\eqref{bou4}.

\begin{thm}
\label{thm3.1}
Let $B(1)\subset \Bbb R^2 $ be the unit ball, and $ a,b, d, e\in C^2(\overline{B(1)} )$ be the coefficients of the system  \eqref{bou1}-\eqref{bou4}. 
Let $ 2< q < +\infty$, and 
\begin{align*}
(u, \theta ) &\in  C([-1,0); W^{2,\, q}(B(1)))\times C([-1,0); W^{2,\, q}(B(1))),
\\
  w  &\in  C([-1,0); W^{1,\, q}(B(1)))
\end{align*}
be a solution to  \eqref{bou1}-\eqref{bou4}. Suppose that  
\begin{equation}
u\in C_w([-1, 0]; L^2(B(1)))
\label{type1-bou-a}
\end{equation}
and  
\begin{equation}
 \intl_{-1}^{0} (-t) \|\nabla  \theta (t)\|_{ L^\infty(B(1))}   dt <+\infty, \quad 
   \intl_{-1}^{0}  \| u(t)\|_{L^\infty (B(1))}dt <+\infty. 
\label{type1-bou}
\end{equation}
Then  $ u, \theta \in C([-1, 0], W^{2,\, q}(B(r)))$ for all $ 0< r< 1$.
\end{thm}
\begin{rem}
As is shown in the next section,  the statement of  Theorem\,\ref{thm1.1}  reduces to the statement of Theorem\,\ref{thm3.1} 
by means of a suitable transformation,  where the condition \eqref{type1-bou-a} is satisfied  by  the global  assumption 
$ v\in L^\infty(-1,0; L^2(\R^{3}))$ of the solution to the 3D Euler equations \eqref{euler}.   
\end{rem}

For the proof of Theorem\,\ref{thm3.1} we will make use of the following Lemma.

\begin{lem}
\label{lem3.3}
Under the assumption of Theorem\,\ref{thm3.1} it holds for every $ 0<r< 1$
\begin{align}
 \intl_{-1}^{0}  \|  w  (t)\|_{ L^\infty(B(r))}  dt < +\infty, 
\label{3.6}
\\
\limsup_{t \to 0^-} (-t) \|  w  (t)\|_{ L^\infty(B(r))} =0,
\label{3.7}
\\
\sup_{t \in (-1,0)}  \| \theta  (t)\|_{ L^\infty(B(r))} < +\infty.
\label{3.7-th}
\end{align}
\end{lem}

{\bf Proof}:  Let $  \frac{1}{2}<r< 1$ be fixed.  Due to \eqref{type1-bou} we may  choose $ -1< t_0<0$ such that 
\begin{equation}
  \intl_{t_0}^{0} \| u(s)\|_{ L^\infty(B(1))}  ds \le 
  \frac{1}{40 a_0} \frac{1-r}{1+r},
\label{3.7a}
\end{equation}
where $ a_0 = \max_{y\in B(1)} | a(y)|$.  Let $ r_{0} = \frac{1+r}{2}$. For $(y,t) \in B(r_0)\times (t_0, 0)$ we define  the coordinate transformation $(y,t)\mapsto  (\widetilde{\varrho}  (y,t), t)$ by
 \[
\widetilde{\varrho}  (y,t)= \varrho(t)y \quad \text {with}\quad \varrho(t):= 1+  20a_0\intl_{t}^{0} \| u(s)\|_{ L^\infty(B(1))}  ds.
\]
One can check easily that 
\begin{align*}
\widetilde{\varrho}  (y,t) &\in  B(1)\quad  \forall\, (y,t) \in B(r_0)\times (t_0, 0), 
\\
\partial _t \widetilde{\varrho}  (y,t) &=-20a_0 \| u(t)\|_{ L^\infty(B(1))} y. 
\end{align*}
Under the above coordinate transform we set
 \begin{align*}
U (y,t) &= u( \widetilde{\varrho}  (y,t), t),
\\
\Theta(y,t)  &= \theta ( \widetilde{\varrho} (y,t), t),
\\
\Omega(y,t) &=  w   ( \widetilde{\varrho}  (y,t), t), \quad  (y,t)\in B(r_0)\times (t_0, 0). 
\end{align*}
In addition, define  
\begin{equation}\label{defw}
W(y ,t) =  \frac{ 20a_0\| u(t)\|_{ L^\infty(B(1))} y  + a( \widetilde{\varrho}  (y,t))U(y,t)}{\varrho (t)},\quad  (y,t) \in B(r_0)\times (t_0, 0). 
\end{equation}

We claim that   
\begin{equation}
 y\cdot W(y,t) \ge \frac{a_0}{8}\| u(t)\|_{ L^\infty(B(1))} \quad  \forall\, (y,t)\in \overline{B(r_0)} 
  \setminus  B(r/2) \times (t_0, 0). 
\label{3.5}
\end{equation}
We now estimate   
\begin{align*}
\varrho(t ) y\cdot W(y,t)&=
20a_0\| u(t)\|_{ L^\infty(B(1))}| y|^2 +   a(\widetilde{\varrho}  (y,t))y\cdot U(y, t)
\\
&\ge \frac{5}{4} a_0\| u(t)\|_{ L^\infty(B(1))} -  a_0| U(y, t)|
\\
&\ge \frac{1}{4} (5a_0\| u(t)\|_{ L^\infty(B(1))} -  4a_0\| u(t)\|_{ L^\infty(B(1))})
\\
&=\frac{a_0}{4} \| u(t)\|_{ L^\infty(B(1))}. 
\end{align*}
From \eqref{3.7a} we find $ \varrho(t ) \leq \frac{3+r}{2+2r} \leq \frac43$,  and therefore \eqref{3.5} follows.

\vspace{0.3cm}
Using the chain rule, we see that \eqref{bou2} turns into the following equations hold in $ B(r_0)\times (t_0,0)$
\begin{equation}
\partial _t \Omega + W\cdot \nabla \Omega =  \frac{b(\widetilde{\varrho}  (y,t))}{ \varrho (t)} \cdot  \nabla  \Theta.  
\label{3.10}
\end{equation}

Given $ (y, s)\in B(r_0)\times (t_0, t)$,  we denote    by  $X(\cdot)=  X(y,s; \cdot ): [t_0, s] \rightarrow  \R^{2}$    the trajectory  such that 
\begin{equation}
\label{traj}
\dot{X}(\tau ) = W(X(\tau ), \tau ),\quad  X(s) = y. 
\end{equation}

 We claim that $ X(\tau ) \in B(r_0)$ for all $ \tau \in [t_0, s]$. 
Otherwise there exists $ \tau _0\in [t_0, s]$ such that $ X(\tau _0)\in \partial B(r_0)\times (t_0,0)$ and $ X(\tau )\in B(r_0)$ 
for all $ \tau \in (\tau_0, s]$.  In particular, the function $\tau \mapsto |X(\tau)|^2$ is non increasing at $\tau =\tau_0$. This gives  
\[
0\geq  \frac{d}{d\tau} |X(\tau)|^2 \big|_{\tau=\tau_0}=2X(\tau _0)\cdot \dot{X}(\tau _0)  = 2X(\tau _0)\cdot  W(X(\tau _0), \tau _0) ,
\]
which contradicts \eqref{3.5}.  On the other hand, by the chain rule \eqref{3.10} gives 
\[
 \frac{d}{d\tau }  \Omega (X(\tau ), \tau ) = \frac{b(\widetilde{\varrho}  (X(\tau ), \tau ), \tau )}{\varrho (\tau )} 
 \nabla  \Theta(X(\tau ), \tau ). 
\]   
Recalling that $ X(s)=y$, integration over $ (s_0, s), t_0 \le s_0 \le t$ yields 
\[
\Omega (y, s) = \Omega (X(s_0), s_0) +  \intl_{s_0}^{s} \frac{b(\widetilde{\varrho}  (X(\tau ), \tau ), \tau)}{ \varrho ( \tau )} 
 \nabla  \Theta(X(\tau ), \tau ) d\tau .
\]
Accordingly,   
\begin{equation}
\| \Omega (s)\|_{ L^\infty(B(r_0))} \le \| \Omega (s_0)\|_{ L^\infty(B(r_0))} 
+ b_0 \intl_{s_0}^{s} \| \nabla   \Theta(\tau )\|_{ L^\infty(B(r_0))} d\tau,
\label{3.11}
\end{equation}
where $ b_0 = \max_{y\in B(1)} | b(y)|$. 
In \eqref{3.11} we insert  $ s_0=t_0$,  integrate both sides  over $ (t_0, t)$. This, together with the integration by part gives
\begin{align*}
& \intl_{t_0}^{t}  \| \Omega (s)\|_{ L^\infty(B(r_0))}  ds
\\
& \le  (t-t_0 ) \| \Omega (t_0)\|_{ L^\infty(B(r_0))} 
- b_0\intl_{t_0}^{t} (-s)'\intl_{s_0}^{s} \|   \nabla   \Theta(\tau )\|_{ L^\infty(B(r_0))} d\tau ds
\\
& \le  (t-t_0) \| \Omega (t_0)\|_{ L^\infty(B(r_0))} 
+b_0\intl_{t_0}^{t} (-s) \|  \nabla  \Theta(s)\|_{ L^\infty(B(r_0))}  ds
\\
& \le  (t-t_0) \|  w   (t_0)\|_{ L^\infty(B(1))} 
+b_0\intl_{t_0}^{t} (-s)  \| \nabla   \theta (s)\|_{ L^\infty(B(1))}  ds.
\end{align*}
Since $\| \Omega (s)\|_{ B(r_0)} \ge \|  w   (s)\|_{ B(r)}$, this proves \eqref{3.6}. To verify \eqref{3.7}, we first note that  \eqref{3.11} multiplied by $ (-s)$ implies   
\begin{align}
(-s)\|  w  (s)\|_{ L^\infty(B(r))} &\le (-s)\| \Omega (s)\|_{ L^\infty(B(r_0))} 
\cr
&\le   (-s)\|  w  (s_0)\|_{ L^\infty(B(r_0))} 
+ b_0\intl_{s_0}^{0} (-\tau )\| \nabla \theta(\tau )\|_{ L^\infty(B(r_0))} d\tau. 
\label{3.12}
\end{align}
Applying $ \limsup$ as $ s \rightarrow 0^-$ to both sides of \eqref{3.12}, we are led to 
\begin{align*}
&\limsup_{ s\to 0^-}(-s)\|  w  (s)\|_{ L^\infty(B(r))} 
\\
&\le \limsup_{ s\to 0^-}(-s)\|  w  (s)\|_{ L^\infty(B(1))} \le b_0\intl_{s_0}^{0} (-\tau )
\| \nabla \theta(\tau )\|_{ L^\infty(B(1))} d\tau <+\infty. 
\end{align*}
On the other hand, the integral on the  right-hand side of this inequality tends to $ 0$ as $ s_0 \rightarrow 0$, and  we  conclude  
\eqref{3.7}.  

\hspace{0.5cm}
To prove \eqref{3.7-th} we argue as above.  In fact,  for the same trajectory as in \eqref{traj}  we deduce from  \eqref{bou1} that
\bb\label{3.12a}
\frac{d}{d \tau } \Theta (X(\tau), \tau ) = 0\quad \text{ for all}\quad  \tau \in [t_0, s].    
\ee
Integration of \eqref{3.12a}  over $ (t_0, s)$ gives  
\[
\Theta (X(s),s )=   \Theta(X(t_0), t_0 )= \theta(\widetilde{\varrho} (X(t_0), t_0), t_0 ).    
\]
Accordingly, for all $r<1$ we have
\[
\| \theta (s) \|_{ L^\infty(B(r))} \le \| \theta (t_0)\|_{ L^\infty(B(1))}. 
\]
This  completes the proof of the lemma. \hfill \Beweisende

\vspace{0.5cm}  
{\bf Proof of Theorem\,\ref{thm3.1}.} 
Applying $ \partial _i$, $ i=1,2$, to \eqref{3.10} we get 
\begin{equation}
\partial _t \partial _i\Omega + W\cdot \nabla \partial _i\Omega =
-\partial _i W\cdot \nabla \Omega +   \frac{b(\widetilde{\varrho}  (y,t))}{\varrho (t)}\cdot  \partial _i\nabla  \Theta  
+ \frac{\partial _i b(\widetilde{\varrho}  (y,t))}{\varrho (t)}  \cdot   \nabla  \Theta.   
\label{3.13}
\end{equation}   

Let $ \frac{1}{2} < r < 1$, and choose  $ -1 < t_0 < 0$ to satisfy \eqref{3.7a}.    Set $ \rho _{ \ast} = \frac{r+ r _0}{2}$, and define
\[
r_m := \rho _{ \ast}- (\rho_{ \ast}-r)^{ m+1} \rho _{ \ast}^{ -m}, \quad  m \in \N\cup \{ 0\}. 
\]
Clearly, 
\[
r_{m+1}- r_{ m} = r \left(1- \frac{r}{\rho_* }\right)^{ m+1},\quad\text{and}\quad r_m\nearrow \rho_*.
\]

Let $ \eta_m \in  C^{\infty} (\R)$ denote a cut off function such that $ 0 \le \eta_m  \le  1 $ 
in $ \R$, $ \eta_m \equiv  1$ on $ (-\infty, r_{ m}]$, $ \eta_m \equiv 0$ in $ (r_{ m+1}, +\infty)$, and 
$ 0 \le  -\eta '_m \le \frac{2}{r_{ m+1}- r_{ m}} = 2 r^{ -1} \Big(\frac{\rho _0+r}{\rho _0-r}\Big)^{m+1}$. 
Let $q>2$ be fixed.  We multiply both sides of \eqref{3.13} by $ \partial _i \Omega | \nabla \Omega |^{ q-2} \phi_m^{ 6q}$,  where 
 $\phi_m (y) = \eta_m (| y|) $, integrate the result over  $ B(r_{ m+1})\times (t_0, t)$, $ t_0 < t < 0$, and sum it over $i=1,2$. Then, applying the integration by parts, we have
 \begin{align}
 & \| \nabla \Omega(t) \phi _m^6 \|_{ q}^{ q}  
 -6q\intl_{t_0}^{t} \intl_{B( r_{ m+1})  \setminus B(r_m)} \frac{W(s)\cdot y }{|y|}  
 | \nabla \Omega(s) |^q \phi _m^{ 6q-1} \eta_m '(| y|) dy ds
 \cr &\qquad  = \| \nabla \Omega(t_0) \phi_m  ^6\|_{ q}^{q} +  I+II+III+IV,
\label{3.14}
 \end{align}
 where we set
  \begin{align*}
 &I+II+III+IV 
 \cr
 &= \intl_{t_0}^{t} \intl_{B(r_{ m+1})}
\nabla \cdot W(s) | \nabla \Omega(s) |^q  \phi _m ^{ 6q} dyds 
\cr
&\qquad -q\intl_{t_0}^{t} \intl_{B(r_{ m+1})}
\nabla W(s): \nabla \Omega(s) \otimes \nabla \Omega(s)  | \nabla \Omega(s) |^{ q-2} \phi_m ^{ 6q} dyds 
\cr
&\qquad + q\intl_{t_0}^{t} \intl_{B(r_{ m+1})}
\frac{b(\widetilde{\varrho}  (y, s))}{\varrho (s)}\nabla ^2 \Theta(s) \cdot \nabla  \Omega (s)  | \nabla \Omega(s) |^{ q-2} \phi_m ^{ 6q} dyds
\cr
&\qquad  + q\intl_{t_0}^{t} \intl_{B(r_{ m+1})}
\frac{\nabla b(\widetilde{\varrho}  (y, s)) }{ \varrho (s)}   :\nabla  \Theta(s) \otimes \nabla  \Omega (s)  | \nabla \Omega(s) |^{ q-2} \phi_m ^{ 6q} dyds
 \end{align*}
respectively,  where and hereafter we use the notation $ M : N= \sum_{i,j=1}^2 M_{ij}N_{ij}$ for the matrices $M, N$. 
Since  $ B(r_{ m+1})  \setminus B(r_{ m}) \subset B(r_0)  \setminus B(r/2) $, \eqref{3.5} ensures that 
  $ W\cdot x >0 $ in $ B(r_{ m+1})  \setminus B(r_{ m})\times (t_0 , 0)$. Furthermore, recalling that $ \eta_m ' \le 0$, we see that 
  the sign of the  integral on the left-hand side of \eqref{3.14}   is non-negative.  Consequently,   it follows that 
 \begin{align}
  \| \nabla \Omega(t) \phi_m  ^6\|_{ q}^q  
  \le  \| \nabla \Omega(t_0) \phi _m^6 \|_{ q}^q  +  I+ II+ III + IV.
\label{3.15}
 \end{align}

 We first  focus  on the estimate of $ \| \nabla^2 U(t) \phi_m  ^6\|_{ q}^q$ in terms of $ \| \nabla \Omega(t) \phi_m  ^6\|_{ q}^q$ in order estimate the left-hand side of \eqref{3.15} from below. Recalling the definition of $ U$ and $ \Omega $, by means of  the chain rule,  the equation \eqref{bou4} becomes 
 \begin{equation}
 \begin{cases}
 \curl U(y,t) =  \varrho (t) \left\{d(\widetilde{\varrho}  (y,t)) \Omega(y,t) +  e(\widetilde{\varrho}  (y,t)) \cdot  U(y,t)\right\}
  \\[0.3cm]
 (y,t)\in B(r_0)\times (t_0,0).  
 \end{cases}
 \label{bou4-trans}
 \end{equation}

In view of \eqref{3.45} with $ u=U(t)$,  $ \psi =\phi _m$ and $ k=6q$ we estimate 
\begin{equation}
\| \nabla^2 U(t) \phi_m  ^6\|_{ q}^q\le  c \| \nabla \curl U(t) \phi _m^{ 6}\|^q_{ q} +  c\| \nabla \phi _m \|_{ \infty}^{ 3q-2} 
\| U(t) \phi _m^{ 2}\|^q_{2}. 
\label{3.46}
\end{equation}
From \eqref{bou4-trans} deducing 
\begin{align*}
\nabla \curl U &=  \varrho (t) d( \widetilde{\varrho}( \cdot) ) \nabla \Omega+  e(\widetilde{\varrho}( \cdot) ) \cdot  \nabla U
\\
&\qquad +\varrho (t) ^2 \nabla d(\widetilde{\varrho}( \cdot)  ) \Omega+  \varrho (t)^2 \nabla e(\cdot )  U, 
\end{align*} 
we estimate in $ B(r_0)\times (t_0, 0)$
\begin{align*}
| \nabla \curl U| \le  c(| \nabla \Omega|+ | \nabla U| + | U|),
\end{align*} 
with a constant $ c>0$ depending only on $ \| d\|_{ \infty}, \| \nabla d\|_{ \infty}, \| e\|_{ \infty}$ and $ \| \nabla d\|_{ \infty}$.
Accordingly, 
\begin{equation}
\| \nabla \curl U (t) \phi _m^{ 6}\|_{ q} \le  c \| \nabla \Omega (t) \phi _m^{ 6}\|_{ q} +c \| \nabla U(t) \phi _m^6\|_{ q} + c \| U(t) \phi _m^{ 6}\|_{q}. 
\label{3.48}
\end{equation}
By means of \eqref{3.44} with $ m= k=6$ we find
\begin{equation}
\| \nabla  U(t) \phi _m^{ 6}\|_{q} \le c \|  U(t) \phi _m^6\|^{ \frac{1}{2}}_{ q} \| \nabla  \curl U(t) \phi _m^6\|^{ \frac{1}{2}}_{ q}+ c\| \nabla \phi _m\|_{ \infty}\| U(t) \phi _m^5\|_{q}. 
\label{3.49}
\end{equation}
Estimating the second term on the right-hand side of \eqref{3.48} using \eqref{3.49} along with Young's inequality,  
 it follows that 
\begin{align}
&\| \nabla \curl U (t) \phi _m^{ 6}\|_{ q}
   \le  c \| \nabla \Omega (t) \phi _m^{ 6}\|_{ q} +  
c(1+\| \nabla \phi _m\|_{ \infty})\| U(t) \phi _m^5\|_{q}. 
\label{3.50}
\end{align}

Furthermore, employing \eqref{3.31} with $ m=6$  and $ k=1$, we get 
 \[
\| U(t) \phi _m^{ 5}\|^q_{q} \le c \| U(t)\phi _m^4\|_{2}^{ \frac{q'}{2}} 
\| \nabla U(t) \phi _m^6\|^{ 1- \frac{q'}{2}}_{ q} + \| \nabla \phi _m\|_{ \infty}^{ \frac{q-2}{q}} 
\| U(t)\phi _m^4\|_{2},
\] 
 where  $q':=\frac{q}{q-1}$. 
This inequality combined with \eqref{3.49} yields 
\begin{align}
&\| U(t) \phi _m^{ 5}\|_{q} 
\cr
&\le c \|  U(t) \phi _m^6\|^{ \frac{1}{2}- \frac{q'}{4}}_{ q}   \| \nabla  \curl U(t) \phi _m^6\|^{ \frac{1}{2}- \frac{q'}{4}}_{ q} \| U(t)\phi _m^4\|_{2}^{ \frac{q'}{2}} 
\cr
&\qquad + c\| \nabla \phi _m\|_{ \infty}^{ 1- \frac{q'}{2}}\| U(t) \phi _m^5\|^{ 1- \frac{q'}{2}}_{q} 
\| U(t) \phi _m^5\|^{  \frac{q'}{2}}_{2}    + c\| \nabla \phi _m\|_{ \infty}^{ \frac{q-2}{q}} 
\| U(t)\phi _m^4\|_{2}. 
\label{3.51}
\end{align}
 Applying Young's inequality, we infer from \eqref{3.51}
   \begin{align}
&\| U(t) \phi _m^{ 5}\|_{q} 
\cr
&\le c    \| \nabla  \curl U(t) \phi _m^6\|^{ \frac{2-q'}{2+q'}} _{ q} 
\| U(t)\phi _m^4\|_{2}^{ \frac{2q'}{2+q'}} 
 +  c\| \nabla \phi _m\|_{ \infty}^{ \frac{q-2}{q}} 
\| U(t)\phi _m^4\|_{2}. 
\label{3.52}
\end{align}
  Finally, inserting \eqref{3.52} into the right-hand side of \eqref{3.50} and using Young's inequality, we obtain    
\begin{align}
&\| \nabla \curl U (t) \phi _m^{ 6}\|_{ q}
   \le  c \| \nabla \Omega (t) \phi _m^{ 6}\|_{ q} +  
c(1+\| \nabla \phi _m\|_{ \infty})^{ \frac{2}{q'}})\| U(t) \phi _m^4\|_{2}. 
\label{3.53}
\end{align}
 Combining \eqref{3.46} and  \eqref{3.53},   noting that 
 $$ (1+\| \nabla \phi _m\|_{ \infty}) \le c (r_{ m-1}- r_m)^{ -1}, \quad
  \| U(t)\|_{ L^2(B(r_{ m+1}))} \le \| u\|_{ L^\infty(-1,0; L^2(B(1)))} <+\infty
  $$
  by the assumption of the theorem, we get 
 \begin{equation}
 \| \nabla^2 U(t) \phi_m  ^6\|_{ q}^q\le c \| \nabla \Omega (t) \phi _m^{ 6}\|^q_{ q} +  c (r_{ m-1}- r_m)^{ - 3q+2}
 \label{3.54}
 \end{equation}
  with a constant $ c>0$ independent of $ t\in (t_0, 0)$ and $ m\in \N$.

We continue to estimate the right-hand side of \eqref{3.15} from above. Calculating 
\begin{align*}
\nabla \cdot W(y,t) &= \frac{40a_0 \| u(t)\|_{ L^\infty(B(1))}}{ \varrho (t)} + \nabla a(\widetilde{\varrho}  (y,t))\cdot U(y,t),
\\ 
\nabla W(y,t)&=\frac{20a_0 \| u(t)\|_{ L^\infty(B(1))}}{ \varrho (t)}I + \nabla a(\widetilde{\varrho}  (y, t)) \otimes  U (y,t)
+ \frac{a(\widetilde{\varrho}  (y, t))}{ \varrho (t)} \nabla U(y,t),
\end{align*}
where $I$ denotes the $2\times2$ unit matrix,  we easily get 
\begin{align*}
I+ II &= (2-q)  \intl_{t_0}^{t} \intl_{B(r_{ m+1})}
\frac{20a_0 \| u\|_{ L^\infty(B(1))}}{ \varrho } | \nabla \Omega |^q  \phi _m ^{ 6q} dyds 
\\
&\quad + \intl_{t_0}^{t} \intl_{B(r_{ m+1})}
\nabla a(\widetilde{\varrho}( \cdot) ) \cdot U | \nabla \Omega |^q  \phi _m ^{ 6q} dyds 
\\
& \quad -q \intl_{t_0}^{t} \intl_{B(r_{ m+1})}
\nabla a(\widetilde{\varrho}( \cdot)  )\cdot \nabla \Omega  U\cdot \nabla \Omega | \nabla \Omega |^{ q-2}  \phi _m ^{ 6q} dyds 
\\
&\quad  +\intl_{t_0}^{t} \intl_{B(r_{ m+1})}
\frac{a(\widetilde{\varrho}( \cdot) )}{ \varrho } \nabla U: \nabla \Omega\otimes \nabla \Omega | \nabla \Omega |^{ q-2}  \phi _m ^{ 6q} dyds. 
\end{align*}
Since the first term of the right-hand side is non-positive, we find 
\[
I+II \le c\intl_{t_0}^{t} \intl_{B(r_{ m+1})}
 (| U| + | \nabla U|)| \nabla \Omega |^q  \phi _m ^{ 6q} dyds, 
\]

where $ c>0$ depends on $ q, \| a\|_{ \infty}$ and $ \| \nabla a\|_{ \infty}$. Furthermore, it is readily seen that 
\begin{align*}
III + IV 
\le c\intl_{t_0}^{t} \intl_{B(r_{ m+1})}
(| \nabla \Theta |+| \nabla ^2 \Theta| )| \nabla \Omega |^{ q-1} \phi_m ^{ 6q} dyds,
\end{align*}
where $ c>0$ depends on $ q, \| b\|_{ \infty}$ and $ \| \nabla b\|_{ \infty}$. Inserting the above estimates of $ I$, 
$ II$, $ III$, and $ IV$ into \eqref{3.15}, we are led to 
\begin{align}
&\| \nabla \Omega(t) \phi_m  ^6\|_{ q}^q  
\cr
&  \le  \| \nabla \Omega(t_0) \phi _m^6 \|_{ q}^q  +c\intl_{t_0}^{t} \intl_{B(r_{ m+1})}
 (| U| + | \nabla U|)| \nabla \Omega|^q  \phi _m ^{ 6q} dyds 
\cr
&\qquad  + c\intl_{t_0}^{t} \intl_{B(r_{ m+1})}
(| \nabla \Theta |+| \nabla ^2 \Theta| )| \nabla \Omega|^{ q-1} \phi_m ^{ 6q} dyds.
\label{3.16}
\end{align} 
For $ 0< \var <1$ and $  t\in (t_0, 0)$ we set
$$ Z_m(t)=Z^{ U}_m(t)+ Z^{ \Theta }_m(t)$$
with
\begin{align*}
Z^{ U}_m(t):=  \| \nabla ^2 U(t) \phi _m ^6\|_{ q}^q,\qquad
Z^{ \Theta }_m(t):=\var ^{ -q}(-t)^{ q}\| \nabla ^2 \Theta (t) \phi _m^6\|_{ q}^q.
\end{align*}
By  Young's inequality we find 
\begin{align*}
| \nabla ^2 \Theta(y,s)|   | \nabla \Omega(y, s) |^{ q-1}&=\var  (-s)^{ -1} \var^{-1}  (-s) | \nabla ^2 \Theta(y, s)|\} \{| \nabla \Omega(y, s) |^{ q-1}
\\
&\le \var (-s)^{ - 1} \Big\{\var ^{ -q}(-s)^{ q}| \nabla ^2 \Theta(y, s)| ^q  + | \nabla \Omega(y, s) |^{ q} \Big\}.
\end{align*}
Therefore,
\bb\label{3.16b}
\intl_{t_0}^{t} \intl_{B(r_{ m+1})}
| \nabla ^2 \Theta|| \nabla \Omega |^{ q-1} \phi_m ^{ 6q} dyds
\leq \var  \intl_{t_0}^{t} (-s)^{ -1} Z_m  (s) ds,
\ee
and
\begin{align}
\intl_{t_0}^{t} \intl_{B(r_{ m+1})}
| \nabla \Theta || \nabla \Omega |^{ q-1} \phi_m ^{ 6q} dyds&= \intl_{t_0}^{t} \intl_{B(r_{ m+1})}
| \nabla \Theta |\phi_m ^{6}| \nabla \Omega |^{ q-1} \phi_m ^{ 6q-6} dyds
\cr
&\leq   \intl_{t_0}^{t}  \| \nabla \Theta (s) \phi _m^{ 6}   \|_q     \|\nabla \Omega (s) \phi_m ^6 \|_{q}^{q-1} ds.
\label{3.16c}
\end{align}
Combining \eqref{3.16} with \eqref{3.16b} and \eqref{3.16c},  we get 
 \begin{align}
 \| \nabla \Omega(t) \phi_m  ^6\|_{ q}^q  &\le  cZ^U_m(t_0)  +c\intl_{t_0}^{t}
 \Big[\| U(s)\|_{ L^\infty(B(r_{ m+1}))} + \| \nabla U(s)\|_{ L^\infty(B(r_{ m+1}))} \Big]
 Z_m^U  (s) ds
\cr
& \qquad + \kappa \var  \intl_{t_0}^{t} (-s)^{ -1} Z_m  (s) ds 
  + c\intl_{t_0}^{t}  \| \nabla \Theta (s) \phi _m^{ 6}\|_{ q}
Z_m^U(s)^{ \frac{q-1}{q}} ds
\label{3.16a}
\end{align} 
with a constant $ \kappa>0$ independent of  $ \var $. 

\hspace{0.5cm}
Now we turn to the estimation of second gradient of $ \Theta $.   First, by the chain rule we derive  from \eqref{bou2} the 
following equation for $ \Theta $ in $ B(r_0)\times (t_0,0)$
\begin{equation}
 \partial _t \Theta   + W\cdot \nabla \Theta  = 0.  
\label{3.17}
\end{equation}

We apply, $ \partial _i\partial _j $, $ i,j=1,2$ to \eqref{3.18},  to obtain 
 \begin{equation}
 \partial _t \partial _i\partial _j \Theta   + W\cdot \nabla \partial _i\partial _j \Theta  = 
 -\partial _iW\cdot \nabla \partial _j \Theta -  \partial _jW\cdot \nabla \partial _i \Theta  - \partial _i\partial _j W\cdot \nabla \Theta .  
\label{3.18}
\end{equation}

We now multiply both sides of \eqref{3.18} by $ \var ^{ -q}(-s)^{ q}\partial _i\partial _j\Theta  | \nabla^2 \Theta  |^{ q-2} \phi_m^{ 6q}$, 
integrate the result over  $ B(r_{ m+1})\times (t_0, t)$, $ t_0 < t < 0$, and sum it over $i=1,2$. Then, applying the integration by parts, we have
 \begin{align*}
 &\var ^{ -q}(-t)^q \| \nabla^2 \Theta (t) \phi _m^6 \|_{ q}^{ q}  
 \cr
 &\qquad \qquad -6q\intl_{t_0}^{t} \intl_{B( r_{ m+1})  \setminus B(r_m)} (-s)^{ q}\frac{W\cdot y }{|y|}  
| \nabla^2\Theta |^q \phi _m^{ 6q-1} \eta_m '(| y|) \var ^{ -q}(-s)^{ q}dyds
 \end{align*}
  \begin{align}
   &\quad =\var ^{ -q}(-t_0)^{ q}  \| \nabla ^2\Theta (t_0) \phi_m  ^6\|_{ q}^{q} + \intl_{t_0}^{t} \intl_{B(r_{ m+1})}
\nabla \cdot W | \nabla^2 \Theta |^q  \phi _m ^{ 6q} \var ^{ -q}(-s)^{ q}dyds 
\cr
&\qquad -2q\intl_{t_0}^{t} \intl_{B(r_{ m+1})}
\nabla W: \nabla^2 \Theta \otimes \nabla^2 \Theta | \nabla^2 \Theta  |^{ q-2} \phi_m ^{ 6q}
\var ^{ -q}(-s)^{ q}dyds 
\cr
&\qquad - q\intl_{t_0}^{t} \intl_{B(r_{ m+1})} \nabla ^2 W \cdot \nabla \Theta :
\nabla ^2 \Theta  | \nabla ^2\Theta |^{ q-2} \phi_m ^{ 6q} \var ^{ -q}(-s)^{ q}dyds
\cr &\quad = \var ^{ -q}(-t_0)^{ q}\| \nabla^2 \Theta (t_0) \phi_m  ^6\|_{ q}^{q} +  I+II+III.
\label{3.19}
 \end{align}
 
 Once more, using the fact that $ y\cdot W(y, s) >0$ for all $(y,s)\in  B(r_0)\times (t_0,0)$, we get from \eqref{3.19} 
 \begin{align}
Z_m^{\Theta }(t)   
  \le Z_m^{\Theta }(t_0)   +  I+II+III.
\label{3.19e}
 \end{align}

 Arguing similarly to the above, we estimate 
\begin{align}
I+II &\le  c\intl_{t_0}^{t} \intl_{B(r_{ m+1})}
 (| U(s)| + | \nabla U(s)|)| \nabla^2 \Theta |^q  \phi _m ^{ 6q}\var ^{ -q} (-s)^{ q}dyds 
\cr
&\le c\intl_{t_0}^{t}
 \Big[\| U(s)\|_{ L^\infty(B(r_{ m+1}))} + \| \nabla U(s)\|_{ L^\infty(B(r_{ m+1}))} \Big]
 Z^\Theta _m (s) ds.
 \label{3.19a}
\end{align}
Furthermore, computing from \eqref{defw}
\begin{align*}
\nabla ^2 W 
= \frac{1}{\varrho} \left\{ \nabla a(\widetilde{\varrho}( \cdot) ) \otimes  \nabla U + \nabla U \otimes \nabla a(\widetilde{\varrho}( \cdot) )+ 
\nabla^2 a(\widetilde{\varrho}( \cdot) ) U + a(\widetilde{\varrho}( \cdot) )\nabla ^2 U\right\},   
\end{align*}
we immediately get 
\[
III \le c\intl_{t_0}^{t} \intl_{B(r_{ m+1})} (| U|+ | \nabla U|+ | \nabla ^2 U|) | \nabla \Theta | 
 | \nabla ^2\Theta|^{ q-1} \phi_m ^{ 6q} \var ^{ -q}(-s)^{ q}dyds
\]
with a constant $ c>0$, depending on $q,  \| a\|_{ \infty}, \| \nabla  a\|_{ \infty}$ and $ \| \nabla ^2 a\|_{ \infty}$. 
Using the fact 
\begin{align}
 | \nabla ^2 U(y, s) || \nabla ^2 \Theta(y, s)|^{q-1}  &=\var ^{-1}  (-s)\var  (-s)^{-1}  | \nabla ^2 Uy,(s) | | \nabla ^2 \Theta(y, s)|^{q-1}\cr
&\le \var  ^{-1} (-s)\Big\{\var^q  (-s) ^{-q} | \nabla ^2 U(y, s)| ^q  +  |  \nabla ^2 \Theta(y, s)|^{ q} \Big\},
\label{young2}
\end{align}
we find 
\begin{align*}
&\intl_{t_0}^{t} \intl_{B(r_{ m+1})}| \nabla ^2 U|| \nabla \Theta|  | \nabla ^2\Theta  |^{ q-1} \phi_m ^{ 6q} \var ^{ -q}(-s)^{ q}dyds\cr
& \leq  
\var ^{ -1}\intl_{t_0}^{t} (-s)  \| \nabla \Theta (s)\|_{ L^\infty(B(r_{ m+1}))} Z_m (s)  ds.
\end{align*}
Similarly, using  H\"{o}ler's inequality  for the other terms of $III$, we estimate
\begin{align}
III &\le c\var ^{ -1}\intl_{t_0}^{t} (-s)\| \nabla \Theta (s)\|_{ L^\infty(B(r_{ m+1}))}  Z_m(s) ds
\cr
& +c\var ^{ -1}\intl_{t_0}^{t} \| U(s)\|_{ L^\infty(B(r_{ m+1}))} (-s)\| \nabla \Theta(s) \phi _m^6\|_{ q} Z_m^\Theta (s) ^{ \frac{q-1}{q}} ds
\cr
& +c\var ^{ -1}\intl_{t_0}^{t}  \| \nabla U(s)\phi _m^3\|_{ L^q(B(r_{ m+1}))} 
(-s)\| \nabla \Theta(s) \phi _m^3\|_{ \infty} Z_m^\Theta (s) ^{ \frac{q-1}{q}} ds.
\label{3.19b}
\end{align}
Inserting the  estimates  \eqref{3.19a} and \eqref{3.19b} into the right-hand side of \eqref{3.19e}, we arrive at 
\begin{align}
&Z^{ \Theta }_m (t) \le Z^{ \Theta }_m (t_0)  
\cr
&+ c \intl_{t_0}^{t}
 \Big[\| U(s)\|_{ L^\infty(B(r_{ m+1}))} + \| \nabla U(s)\|_{ L^\infty(B(r_{ m+1}))} + (-s)\| \nabla \Theta (s)\|_{ L^\infty(B(r_{ m+1}))}\Big]  Z_m (s) ds
 \cr
& +c\intl_{t_0}^{t} \| U(s)\|_{ L^\infty(B(r_{ m+1}))} (-s)\| \nabla \Theta(s) \phi _m^6\|_{ q} Z_m^\Theta (s) ^{ \frac{q-1}{q}} ds
\cr
& +c\intl_{t_0}^{t}  \| \nabla U(s)\phi _m^3\|_{ L^q(B(r_{ m+1}))} 
(-s)\| \nabla \Theta(s) \phi _m^3\|_{ \infty} Z_m^\Theta (s) ^{ \frac{q-1}{q}} ds
\label{3.21}
\end{align} 
with a constant $c$ depending on $ \var $. 
Taking the sum of \eqref{3.16a} and \eqref{3.21}, and taking into account \eqref{3.54},  we obtain   
 \begin{align}
Z_m(t) &\le cZ_m (t_0)  + c (r_{ m+1}- r_m)^{- 3q+2}
\cr
&+ c \intl_{t_0}^{t}
 \Big[\| U(s)\|_{ L^\infty(B(r_{ m+1}))} + \| \nabla U(s)\|_{ L^\infty(B(r_{ m+1}))} + (-s)\| \nabla \Theta (s)\|_{ L^\infty(B(r_{ m+1}))}
 \Big]  Z _m (s) ds
 \cr
 &  + \kappa\var \intl_{t_0}^{t} (-s)^{ -1} Z _m (s) ds+ J_1+J_2+J_3,
 \label{3.22}
 \end{align}   where we set
 \begin{align*}
 J_1&= c\intl_{t_0}^{t} \| U(s)\|_{ L^\infty(B(r_{ m+1}))} (-s)\| \nabla \Theta(s) \phi _m^6\|_{ q} Z_m^\Theta (s) ^{ \frac{q-1}{q}} ds,\\
J_2&=c\intl_{t_0}^{t}  \| \nabla \Theta (s)\phi _m^6\|_{ L^q(B(r_{ m+1}))} 
 Z_m^U (s) ^{ \frac{q-1}{q}} ds,
 \end{align*} 
 and 
$$
J_3= c\intl_{t_0}^{t}  \| \nabla U(s)\phi _m^3\|_{ L^q(B(r_{ m+1}))} 
(-s) \| \nabla \Theta (s) \phi _m^{ 3}\|_{ q}
Z_m^\Theta(s)^{ \frac{q-1}{q}} ds.
$$ 
As we shall show below,  the  integrals $ J_1, J_2$ and $ J_3$ are lower order terms dominated by the other terms in the right-hand side of \eqref{3.22}. To see this we 
first estimate the term  $ \| \nabla \Theta(s) \phi _m^3\|_{ q}$ as follows. 
Observing \eqref{3.7-th} (cf. Lemma\,\ref{lem3.3}), we easily verify that 
\begin{equation}
\Theta _\infty:= \sup_{ t\in (t_0, 0)}\| \Theta(t) \|_{ L^\infty(B(r_0))} < +\infty.  
\label{3.23}
\end{equation}
Integrating by parts, we calculate 
\begin{align}
&\| \nabla \Theta(s) \phi _m^3\|_{ q}^q= \intl_{B(r_{ m+1})} \nabla \Theta (y, s) \cdot \nabla \Theta (y, s)| \nabla \Theta (y, s)|^{ q-2} \phi _m^{ 3q}(y) dy
\cr
&=- \intl_{B(r_{ m+1})}  \Theta (y,s) \Delta  \Theta (y, s) \phi _m^{ 6}(y)| \nabla \Theta (y, s)|^{ q-2} \phi _m^{ 3(q-2)} (y)dy 
\cr
&\,  - (q-2)\sum_{i=1} ^2 \intl_{B(r_{ m+1})}  \Theta (y, s)  \partial _i \Theta (y, s) \partial _i \nabla \Theta(y, s)
\phi _m^{ 6}(y)\cdot \nabla \Theta (y, s)  | \nabla \Theta (y, s)|^{ q-4} \phi _m^{ 3(q-2)}(y) dy
\cr
& \qquad - 6q \intl_{B(r_{ m+1})}  \Theta (y, s)  \nabla \Theta (y, s)\cdot \nabla \phi _m(y)   | \nabla \Theta (y, s)|^{ q-2} \phi _m^{ 3q-1}(y) dy
:=A+B+C.
\label{3.24} 
\end{align} 
Thanks to \eqref{3.23} along with H\"older's inequality and Young's inequality, we immediately get 
\begin{align*}
A+B  &\le c\Theta _\infty  \intl_{B(r_{ m+1})} | \nabla^2 \Theta (y, s)| \phi _m^{6} | \nabla   \Theta (y, s)|^{ q-2} \phi ^{ 3(q-2)}_m(y) dy 
\\
&\le  c\Theta _\infty   \| \nabla ^2 \Theta (s) \phi _m^{ 6}\|_q  \| \nabla \Theta(s) \phi _m^3\|_{ q}^{ q-2}
\\
&\le c\Theta _\infty^{ \frac{q}{2}}   \| \nabla ^2 \Theta (s) \phi _m^{ 6}\|_q^{ \frac{q}{2}} + \frac{1}{4}\| \nabla \Theta(s) \phi _m^3\|^q_{ q}. 
\end{align*} 
Similarly, by means of H\"older's inequality and Young's inequality,  we infer
\begin{align*}
C &\le   6q \intl_{B(r_{ m+1})} |\Theta (y, s)|  |\nabla \phi _m (y)|  | \nabla \Theta (y, s)|^{ q-1} \phi _m^{ 3q-3} (y) \phi_m ^2 (y) dy
\\
 &\le c \Theta _\infty  \| \nabla \phi _m\|_{ \infty} \| \nabla \Theta(s) \phi _m^3\|^{ q-1}_{ q}\le c \Theta _\infty  (r_{ m+1}- r_m)^{ -1} 
 \| \nabla \Theta(s) \phi _m^3\|^{ q-1}_{ q}\\
 &\le   c \Theta _\infty^{ q} (r_{ m+1}- r_m)^{ -q}  + \frac{1}{4} \| \nabla \Theta(s) \phi _m^3\|^q_{ q}. 
\end{align*} 
Inserting the above estimates of $ A, B$ and $ C$ into the right-hand side of \eqref{3.24}, and absorbing  $ \frac12 \| \nabla \Theta(s) \phi _m^3\|^q_{ q}$ into the left-hand side, we deduce that 
\begin{align}
\| \nabla \Theta(s) \phi _m^3\|_{ q} &\le c \Theta _\infty^{ \frac{1}{2}} \| \nabla^2 \Theta(s) \phi _m^6\|^{ \frac{1}{2}}_{ q}  + 
 c \Theta _\infty(r_{ m+1}- r_m)^{ -1} 
\cr
& = c \Theta _\infty^{ \frac{1}{2}} \var ^{ - \frac{1}{2}}(-s)^{ -\frac{1}{2}} Z_m^{ \Theta }(s)^{ \frac{1}{2q}}+ c \Theta _\infty(r_{ m+1}- r_m)^{ -1} \cr
& \leq c (-s)^{ -\frac{1}{2}} Z_m^{ \Theta }(s)^{ \frac{1}{2q}}+ c (r_{ m+1}- r_m)^{ -1}.
\label{3.25}
\end{align}
Similarly, following the above  calculations    \eqref{3.23}-\eqref{3.25}, we also have
\begin{align}
\| \nabla U(s) \phi _m^3\|_{ q} &\le c \|U(s)\|_{ L^\infty(B(r_{ m+1}))} ^{ \frac{1}{2}} \| \nabla^2 U(s) \phi _m^6\|^{ \frac{1}{2}}_{ q}  + 
 c \|U(s)\| _{ L^\infty(B(r_{ m+1}))} (r_{ m+1}- r_m)^{ -1} 
\cr
& = c \|U(s)\| _{ L^\infty(B(r_{ m+1}))} ^{ \frac{1}{2}} Z_m^{ U}(s)^{ \frac{1}{2q}}+ c \|U(s)\| _{ L^\infty(B(r_{ m+1}))} (r_{ m+1}- r_m)^{ -1}.
\label{3.25a}
\end{align}
Obviously, 
$$  \| \nabla \Theta(s) \phi _m^6\|_{ q} \leq \| \nabla \Theta(s) \phi _m^3\|_{ q}, \quad \text{and}\quad    \| \nabla U(s) \phi _m^6\|_{ q} \leq \| \nabla U(s) \phi _m^3\|_{ q}. $$
Therefore, using \eqref{3.25}, and then applying Young's inequality, we find 
\begin{align*}
J_{ 1} &\le c \intl_{t_0}^{t} \| U(s)\|_{ L^\infty(B(r_{ m+1}))} \Big[ (-s)^{  \frac{1}{2} } (Z^{ \Theta }_m(s))^{ \frac{2q-1}{2q}} + (r_{ m+1}-r_{ m} )^{ -1}(Z^\Theta_m(s))^{ \frac{q-1}{q}}\Big]ds
\\
&\le c  \intl_{t_0}^{t} \| U(s)\|_{ L^\infty(B(r_{ m+1}))}   Z^\Theta_m(s)  ds 
 + c(r_{ m+1}-r_m)^{ -q}  \intl_{t_0}^{t} \| U(s)\|_{ L^\infty(B(r_m))}    ds
 \\
 &\qquad +  c  \intl_{t_0}^{t} \| U(s)\|_{ L^\infty(B(r_{ m+1}))}  (-s)^q ds.
 \end{align*} 
Noting that $ \| U(s)\|_{ L^\infty(B(r_{ m+1}))} \le \| u(s)\|_{ L^\infty(B(1))}$, and observing \eqref{type1-bou}$ _2$, we find 
\begin{align*}
J_{ 1} &\le  c  \intl_{t_0}^{t} \| u(s)\|_{ L^\infty(B(1))} Z_m(s) ds   
 + c(r_{ m+1}-r_m)^{ -q} \intl_{t_0}^{t} \| u(s)\|_{ L^\infty(B(1))}  ds\\
 & \qquad+ c\intl_{t_0}^{t} \| u(s)\|_{ L^\infty(B(1))} ds\\
&\le  c  \intl_{t_0}^{t} \| u(s)\|_{ L^\infty(B(1))} Z_m(s) ds   
 + c(r_{ m+1}-r_m)^{ -q}. 
\end{align*}
Using \eqref{3.25} together with Young's inequality,  we easily get 
\begin{align*}
J_2 & \le c \intl_{t_0}^{t} \Big[  (-s)^{  \frac{1}{2} } (Z^{ \Theta }_m(s))^{ \frac{1}{2q}} (Z^U_m(s))^{ \frac{q-1}{q}} 
+ (r_{ m+1}-r_{ m} )^{ -1}(Z^U_m(s))^{ \frac{q-1}{q}}\Big]ds
\\
&\le c  \intl_{t_0}^{t} \Big[  (-s)^{\frac12}\Big(  Z_m(s)\Big)^{ \frac{2q-1}{2q}} + c (r_{ m+1}- r_m)^{ -1}  
Z_m(s)^{ \frac{q-1}{q}}\Big] ds 
\\
& \le c  \intl_{t_0}^{t} Z_m(s)  ds +  c  \intl_{t_0}^{t}  (-s) ^q ds +  c (r_{ m+1}- r_m)^{ -q} \\
&\le c  \intl_{t_0}^{t} Z_m(s)  ds  + c (r_{ m+1}- r_m)^{ -q}.
\end{align*}

Using \eqref{3.25}, \eqref{3.25a} and Young's inequality, we easily estimate
\begin{align*}
J_3 &\le  c  \intl_{t_0}^{t}  (-s)^{\frac12} \|U(s)\|_{ L^\infty(B(r_{ m+1}))} ^{\frac12}(Z_m ^U (s))^{\frac{1}{2q}} (Z_m ^\Theta (s)) ^{\frac{2q-1}{2q}}  ds\\
     &\qquad+c  \intl_{t_0}^{t}   (-s)^{\frac12}  \|U(s)\|_{ L^\infty(B(r_{ m+1}))} (Z_m ^\Theta (s)) ^{\frac{2q-1}{2q}}   (r_{ m+1} - r_m)^{ -1} ds \\
 &\qquad+ c  \intl_{t_0}^{t}   (-s) \|U(s)\|_{ L^\infty(B(r_{ m+1}))}  ^{\frac12} (Z_m ^\Theta (s)) ^{\frac{q-1}{q}} (Z_m ^U (s))^{\frac{1}{2q}}   (r_{ m+1} - r_m)^{ -1} ds \\
  &\qquad+c  \intl_{t_0}^{t}   (-s) \|U(s)\|_{ L^\infty(B(r_{ m+1}))}  (Z_m ^\Theta (s)) ^{\frac{q-1}{q}}  (r_{ m+1} - r_m)^{ -2} \\
  &\le  c\int_{t_0} ^t \|U(s)\|_{ L^\infty(B(r_{ m+1}))} Z_m(s) ds +  c(r_{ m+1} - r_m)^{ -2}  \\
  &\le  c\int_{t_0} ^t \|u (s)\|_{ L^\infty(B(1))} Z_m(s) ds +  c(r_{ m+1} - r_m)^{ -2}.
  \end{align*}
Inserting the estimates of $ J_1, J_2$ and $ J_3$ into the right-hand side of \eqref{3.22}, 
noting that $ 2q < 3q-2$ for $ q>2$, 
and observing \eqref{type1-bou}$ _1$, 
we arrive at  
  \begin{align}
& Z_m (t) \le c Z_m (t_0)  
\cr
&+ c \intl_{t_0}^{t}
 \Big[ \| \nabla U(s)\|_{ L^\infty(B(r_{ m+1}))}+  (-s)\| \nabla \Theta (s)\|_{ L^\infty(B(r_{ m+1}))}+
 \| u(s)\|_{ L^\infty(B(1))}\Big]  Z _m (s) ds
 \cr
 &\qquad \qquad  + \kappa \var  \intl_{t_0}^{t}
  (-s)^{ -1}  Z _m (s) ds+ c(r_{ m+1} - r_m)^{ -3q+2},
  \label{3.55}
\end{align} 
where $ c$ stands for a  positive constants, depending on $ q, \var, \alpha, C_0,  u, \theta $ but not on $ t$ and $ m$, 
while the constant $ \kappa $  depends on $ q, \alpha, C_0,  u, \theta $ but not on  $\var ,  t$ and $ m$.   

By the similar argument to  \cite[cf. (5.14)]{cha5}, observing \eqref{bou4-trans},  we get 
\begin{align}
&\| \nabla U(s) \|_{ L^\infty(B(r_{ m+1}))} 
\cr
&\le c\Big\{1+ |\curl U (s)|_{ BMO(B(\rho_{ \ast}))} \Big\} \log (e+
Z_{ m+1}(s)) 
\cr
& \le   c\Big\{1+ \| U(s)\|_{ L^\infty(B(\rho _{ \ast}))}+
 \|  \Omega (s)\|_{L^\infty(B(\rho_{ \ast}))} \Big\} \log (e+
Z_{ m+1}(s)).  
\label{3.56}
\end{align}   
Combining \eqref{3.55} and \eqref{3.56}, and noting  that $ (r_{ m+1} - r_{ m})^{ -1} \le c \Big(\frac{r_0+r}{r_0-r}\Big)^m$, we   are led to 
  \begin{align}
& e+Z_m (t) \le c_0 +  d^m + \intl_{t_0}^{t} (\alpha  (s) \log(e+Z_{ m+1}(s))  + f(s)) Z _m (s) ds,
   \label{3.57}
\end{align}   
where 
\begin{align*}
\alpha  (s) &=  c \Big(1+  \| u(s)\|_{ L^\infty(B(1))}+
 \|  \Omega  (s)\|_{L^\infty(B(\rho_{ \ast}))}\Big),
\\
f(s) &= c\Big[(-s)\| \nabla \theta  (s) \|_{ L^\infty(B(1))}+ \| u(s)\|_{ L^\infty(B(1))}\Big]+ \kappa \var (-s)^{ -1},
\end{align*}
\begin{align*}
c_0 &=c (Z_{ m}(t_0)+e) \le  c \Big(1+ \| \nabla ^2 u(t_0)\|^q_{ L^q(B(1))} + c\| \nabla ^2 \theta (t_0)\|^q_{ L^q(B(1))}\Big), 
\\
d &= c\Big(\frac{r_0+r }{r_0-r}\Big)^{ 3q-2}. 
\end{align*}
We now define  
\[
\begin{cases}
Y_m (t):= c_0+  d^m + \intl_{t_0}^{t} \left\{ \alpha(s) \log(e+Z_{ m+1}(s))  + f(s)\right\} Z _m (s) ds,
\\[0.3cm]
 t\in (t_0, 0),\quad m\in \N. 
\end{cases}
\]
Thus, \eqref{3.57} gives $ Z_m (t) \le Y_m(t)$, and therefore
\begin{align}
Y_m' (t) &= \left\{ \alpha(t) \log(e+Z_{ m+1}(s))  + f(s)\right\} Z _m (s) 
\cr
&\le\left\{ \alpha(t) \log(e+Y_{ m+1}(s))  + f(s)\right\} Y _m (s). 
\label{3.59}
\end{align}
Dividing both sides of \eqref{3.59} by $ e+ Y_m(s)$ and  setting 
\[
\beta _m(s) := \log(e+Y_{ m}(s)),\quad  s\in [t_0, 0), 
\] 
we deduce from \eqref{3.59} the following recursive  differential inequality 
\begin{equation}
\beta '_m(s) \le \alpha(t) \beta _{ m+1}(s)  + f(s). 
\label{3.60}
\end{equation} 
We now fix $ t\in (t_0,0)$, and  integrate \eqref{3.60} over $ (t_0, t)$. This  yields 
\begin{align}
\beta _m(t) &\le X_m(t_0)  +  \intl_{t_0}^{t}    \alpha(s) \beta _{ m+1}(s)  ds +  \intl_{t_0}^{t}   f(s) ds
\cr
&= \log(e+ Y_{ m}(t_0))  +  \intl_{t_0}^{t}    \alpha(s) \beta _{ m+1}(s)  ds +  \intl_{t_0}^{t}   f(s) ds
\cr
&= \log(e+ c_0 + d^m)  +  \intl_{t_0}^{t}    \alpha(s) \beta _{ m+1}(s)  ds +  \intl_{t_0}^{t}   f(s) ds
\cr
& \le  m \log d +  g(t) +  \intl_{t_0}^{t}    \alpha(s) \beta _{ m+1}(s)  ds ,
\label{3.61}
\end{align}
where we set
\[
g(\tau ) = 1 + \log (e+c_0) +\intl_{t_0}^{\tau }   f(s) ds,\quad  \tau \in [t_0,0). 
\]

In order to apply   Lemma\,\ref{lemA.4} we still need to check if  \eqref{A.5aa} is fulfilled.   By the assumption of the theorem we have 
\[
\max_{ \tau \in [t_0, t]}Z_m (\tau ) \le \| \nabla ^2 u\|_{ L^\infty(-1, t; L^q(B(1)))} + \var ^{ -1}(-t)\| \nabla ^2 \theta  
\|_{L^\infty(-1, t;  L^q(B(1)))} =: c_1(t)<+\infty
\]
for all $ m\in \N$. Accordingly, for $ \tau \in [t_0, t]$
\[
Y_m(\tau ) \le c_0 +  d^m + \intl_{t_0}^{\tau } (\alpha(s) \log(e+c_1)  + f(s)) c_1 ds \le c_2 + d^{ m},
\]
for some constant $ c_{ 2}>0$ depending on $ t$ but independent on $ \tau $. Thus, for all $ \tau \in [t_0, t]$
\[
\beta_m(\tau ) = \log(e+ Y_m(\tau )) \le \log(e+ c_2+d^m) \le 1+ \log (e+ c_{ 2}) + m \log d. 
\]
Clearly, the condition \eqref{A.5aa} of Lemma\,\ref{lemA.4} is satisfied. 
Indeed,
\[
\max_{ \tau \in [t_0, t]} \beta _m(\tau ) \le  1+ \log (e+ c_{ 2}) + m \log d \le K^m \quad  \forall\,m\in \N. 
\]
with $ K=1+ \log (e+ c_{ 2}) +  \log d$. 
 We are now in a position to apply Lemma\,\ref{lemA.4} with $ C= \log d$,  which shows that 
 \begin{equation}
\beta _0(t) \le  g(t) +\log d  \intl_{t_0}^{t} \alpha(s)  ds  e^{\intl_{t_0}^{t} \alpha(\tau )   d\tau }+ \intl_{t_0}^{t} \alpha(s)  g(s) e^{\intl_{s}^{t} \alpha(\tau )   d\tau } ds.
\label{3.60a}
\end{equation}

\hspace{0.5cm}
Next, applying integration by parts, and noting that $ g'(t)=f(t)$ and  $ g(t_0)=1+ \log(e+ c_2)$, we deduce from \eqref{3.60a}
\begin{align}
\beta _0(t) \le  g(t) +\log d  \intl_{t_0}^{t} \alpha(s)  ds  e^{\intl_{t_0}^{t} \alpha(\tau )   d\tau }+ 
\intl_{t_0}^{t}   f(s) e^{\intl_{s}^{t} \alpha(\tau )   d\tau } ds + c e^{\intl_{t_0}^{t} \alpha(\tau ) d\tau}.
\label{3.61}
\end{align}

Setting 
\[
r_1 = \frac{3r+1}{4},
\]
owing to  \eqref{3.7a}, we see that $ \varrho (y, t) \in B(r_1)$ for all $ t\in (t_0, 0)$, and 
 thanks to Lemma\,\ref{lem3.3} we get 
 \[
 \intl_{t_0}^{0} \| \Omega(t )\|_{ L^\infty(B(r_{ \ast}))}    dt \le  
 \intl_{t_0}^{0} \|  w (t)\|_{ L^\infty(B(r_{1}))}    dt< +\infty.  
\]
This shows that 
\begin{align*}
 a _0 &:=\intl_{t_0}^{0} \alpha(t)  d t
\\
&=   c(-t_0)+ 
  \intl_{t_0}^{0}  \Big[ \| u(t)\|_{ L^\infty(B(1))}+ \|  \Omega  (t)\|_{L^\infty(B(\rho_{ \ast}))} \Big]dt < +\infty. 
\end{align*}
Therefore,  choosing $|t_0|$ small enough, we may assume that $ 
a _0 \le 1$, and recalling the definition of $ g$, we obtain from \eqref{3.61} 
 \begin{align*}
 \beta _0 (t) &\le c + (1+ e)  \intl_{t_0}^{t} f(s)  ds + (\log d)  e + c e
 \cr
& \le c + (1+ e) \intl_{t_0}^{t} f(s) d s.
 \end{align*}
From  the definition of $ f$, observing \eqref{type1-bou}$ _1$, we get a constant $ c>0$ such that 
\[
 \beta _0 (t) \le c -  (1+ e)  \kappa  \var (\log (-t) - \log (-t_0))
 \le c + \log (-t)^{ - (1+ e)  \kappa  \var}.     
\]
Consequently,
\[
e+ Y_0(t) \le  c (-t)^{ - (1+ e)  \kappa  \var}\quad \forall\,t\in [t_0, 0). 
\]
Since $ Z_0 \le Y_0$, it follows that 
\begin{equation}
\| \nabla ^2 u(t)\|_{ L^q(B(r))}^q +(- t)^{ -q} \| \nabla ^2 \theta (t)\|_{ L^q(B(r))}^q \le  Y_0(t) \le  
c (-t)^{ - (1+ e)  \kappa \var}\quad \forall\,t\in [t_0, 0). 
\label{3.63}
\end{equation}
Choosing $ \var >0$ so that $ (1+ e)  \kappa \var \le \frac{1}{2}$, we get from \eqref{3.63}
\begin{equation}
\int_{-1} ^0 \Big(\| \nabla ^2 u(t)\|_{ L^q(B(r))}^q + 
(-t)^q\| \nabla ^2 \theta (t)\|_{ L^q(B(r))} ^q\Big) dt < +\infty.
\label{3.64}
\end{equation} 
By  Sobolev's embedding theorem we get from \eqref{3.64}   and \eqref{type1-bou}$ _2$ together with \eqref{3.7-th} 
for all $ 0< r< 1$
\begin{equation}
 \intl_{-1}^{0} \| \nabla u(t)\|_{ L^\infty(B(r))}dt < +\infty.
\label{integr}
\end{equation} 
 
To complete  the proof of the theorem, we first show that $ \nabla \theta \in L^\infty(B(r)\times (-1,0)) $ for all $ 0<r<1$. 
We apply, $ \partial _i$, $ i\in \{ 1,2\}$, to \eqref{3.17},   which gives 
 \begin{equation}
 \partial _t \partial _i\Theta   + W\cdot \nabla \partial _i\Theta  = 
 -\partial _iW\cdot \nabla \Theta \quad  \text{ in}\quad  B(r_0)\times (t_0, 0).
\label{3.65}
\end{equation}
Multiplying both sides of \eqref{3.65} by   $\partial _i\Theta $, and taking the sum from $ i=1$ to $ 2$, we arrive at 
 \begin{equation}
 \partial _t | \nabla \Theta|^2   + W\cdot \nabla | \nabla \Theta|^2  = 
 -2\nabla W:  \nabla \Theta \otimes  \nabla \Theta.
\label{3.65a}
\end{equation}
Let  $ (y, s)\in B(r_0)\times (t_0, t)$.  We denote    by  $X(\cdot)=  X(y,s; \cdot ): [t_0, s] \rightarrow  \R^{2}$    the trajectory  such that 
\begin{equation}
\label{traj1}
\dot{X}(\tau ) = W(X(\tau ), \tau ),\quad  X(s) = y. 
\end{equation}
As we have seen in the proof of Lemma\,\ref{lem3.3},  $ X(\tau ) \in B(r_0)$ for all $ \tau \in [t_0, s]$. In addition from \eqref{3.65a} together with \eqref{traj1}  we deduce that   
\begin{equation}
\frac{d}{d\tau } | \nabla \Theta (X(\tau ), \tau )|^2 = -2\nabla W(X(\tau ), \tau ):  \nabla \Theta(X(\tau ), \tau ) \otimes  \nabla \Theta(X(\tau ), \tau ).
\label{3.66}
\end{equation}
Integrating both sides of \eqref{3.66} over $ (t, s)$ for some $ t\in (t_0, s)$, we obtain 
\begin{align*}
& | \nabla \Theta (X(t ), t)|^2 
\\
& = | \nabla \Theta (X(t_0), t_0)|^2 - 2  \intl_{t_0}^{t}  
\nabla W(X(\tau ), \tau ):  \nabla \Theta(X(\tau ), \tau ) \otimes  \nabla \Theta(X(\tau ), \tau ) d\tau 
\\
& \le| \nabla \Theta (X(t_0), t_0)|^2 + 2  \intl_{t_0}^{t}  
\| \nabla W( \tau )\|_{ L^{ \infty}(B(r_0))} | \nabla \Theta(X(\tau ), \tau )|^2  d\tau.  
\end{align*}
Noting that thanks to \eqref{integr} and \eqref{type1-bou}$ _2$  
 \begin{equation}
 \label{3.66a}
  \int_{t_0} ^ 0  \| \nabla W(\tau )\|_{ L^{ \infty}(B(r_0))} d\tau <+\infty, 
  \end{equation}
 and we are in a position to apply Gronwall's Lemma, which shows that 
\begin{align*}
| \nabla \Theta (y,s)|&\le | \nabla \Theta (X(t_0), t_0)|e^{  \intl_{t_0}^{0}  \| \nabla W( \cdot )\|_{ L^{ \infty}(B(r_0))} d\tau }\\
&\leq \| \nabla \Theta (t_0)\|_{ L^\infty(B(r_0))} e^{  \intl_{t_0}^{0}  \| \nabla W( \cdot )\|_{ L^{ \infty}(B(r_0))} d\tau } <+\infty. 
\end{align*}
Whence,
\begin{equation}
\| \nabla \theta \|_{ L^\infty(B(r)\times (-1, 0))} <+\infty. 
\label{3.67}
\end{equation}

Next, multiplying both sides of \eqref{3.13} by  $ \partial _i \Omega  |\nabla \Omega  |^{ q-2}$, $q>3$,  and taking the sum over $ i=1,2$, we get 
\begin{align}
&\frac{1}{q}\partial _t | \nabla \Omega|^q + \frac{1}{q}W\cdot \nabla | \nabla \Omega|^q 
\cr
&= -\nabla  W:\nabla \Omega \otimes \nabla  \Omega | \nabla \Omega  |^{ q-2} +   
\frac{b(\widetilde{\varrho}  (y,t))}{\varrho (t)} \otimes \nabla \Omega:\nabla^2  \Theta   | \nabla \Omega  |^{ q-2}
\cr
&\qquad \qquad + \frac{\nabla  b(\widetilde{\varrho}  (y,t))}{\varrho (t)}  : \nabla  \Theta \otimes \nabla  \Omega | \nabla \Omega |^{ q-2}
\cr
& \le | \nabla  W|| \nabla \Omega|^{ q} + c_1
| \nabla^2  \Theta|   | \nabla \Omega  |^{ q-1} + c_2| \nabla  \Theta|   | \nabla \Omega  |^{ q-1},
\label{3.69}
\end{align}
where $ c_1, c_2$ are constants, depending on $ b$. Integrating both sides of \eqref{3.69} over $ B(r_0)$ and applying integration by parts, 
and H\"older's inequality,   
we are led to 
\begin{align*}
\hspace{-1.in}&\frac{1}{q} \partial _t \| \nabla \Omega(t) \|_{ L^q(B(r_0))}^q + \frac{1}{q}  \intl_{\partial B(r_0)} \frac{y}{r_0}\cdot  W(y,t) 
 | \nabla \Omega(y,t) |^q dS 
 \\
&\quad  \le   \Big(\frac{1}{q}+1\Big)\| \nabla  W(t)\|_{ L^\infty(B(r_0))} \| \nabla \Omega(t)\|^{ q}_{ L^q(B(r_0))} + c_1
\| \nabla^2 \Theta (t)\|_{ L^q(B(r_0))}   \| \nabla \Omega(t)\|^{ q-1}_{ L^q(B(r_0))} 
\\
&\qquad \qquad + c_2\| \nabla \Theta (t)\|_{ L^q(B(r_0))}   \| \nabla \Omega(t)\|^{ q-1}_{ L^q(B(r_0))}.
\end{align*}

First note that, according to \eqref{3.5},  the second term on the left-hand side  is  nonnegative. Furthermore,  thanks to \eqref{3.64} and \eqref{3.67}  the second and the third term on the right-hand side belong to $ L^1(t_0, 0)$, and the function $ \| \nabla W( \cdot )\|_{ L^{ \infty}(B(r_0))}$  
belongs to $ L^1(t_0, 0)$(see \eqref{3.66a}), we are in a position to apply Gronwall's Lemma, which shows that 
\[
\sup_{t\in (t_0, 0) }  \| \nabla \Omega(t) \|_{ L^q(B(r_0))} < +\infty. 
\] 
This together with  $ u\in C([-1, t_0]; W^{2,\, q}(B(1)))$, which is  an assumption of the theorem,  and \eqref{bou4}  shows that 
\[
 \curl u \in  L^\infty(-1, 0; W^{1,q}(B(r)))\quad  \forall\, 0< r< 1.
\]
Employing Lemma\,\ref{lem3.9}, together with the assumption   $ u\in C_w([-1,0]; L^2(B(1)))$, we get   
\begin{equation}\label{3.69a}
 u \in L^\infty(-1, 0; W^{2,q}(B(r)) \quad \forall\,0< r< 1.  
 \end{equation}
Next, we claim 
\begin{equation}\label{3.69b}
u\in C(\overline{B(r)}  \times [-1, 0]) \quad  \forall\,0< r< 1.
 \end{equation}

To see this, let $ \{ (x_k, t_k)\}$ be a sequence in $\overline{B(r)} \times  (-1,0)$, $ 0<r<1$, which converges to $ (x_0, 0)\in \overline{B(r)}\times \{ 0\}$ as $ k \rightarrow \infty$. Since $ W^{1,\, q}(B(r))$ 
is compactly embedded into $ C(\overline{B(r)} )$, we get from \eqref{3.69a} and the assumption   $ u\in C_w([-1,0]; L^2(B(1)))$ 
that $ u(t_k) \rightarrow u(0)$ uniformly in $ \overline{B(r)} $ as $ k \rightarrow +\infty$.  On the other hand, 
by virtue of Sobolev's embedding theorem we see that $ W^{1,\, q}(B(r))$  is continuously embedded into 
$ C^{\gamma }(\overline{B(r)} )$, for $ \gamma = 1- \frac{2}{q}$.  Hence, using triangle inequality, we find 
\begin{align*}
| u(x_k, t_k)- u(x_0, 0)| &\le | u(x_k, t_k)- u(x_0,  t_k)| + |u(x_0, t_k)- u(x_0, 0)|
\\
&\le c| x_k- x_0|^{ \gamma }\| u(t_k)\|_{ W^{1,\, q}(B(r))}+ \| u(t_k)-u(0)\|_{ L^\infty(B(r))}. 
\end{align*}
Clearly, as $ k \rightarrow +\infty$,  the second term on the right-hand side tends to zero, while the first term 
tends to zero thanks  \eqref{3.69a}. Whence, \eqref{3.69b}.

 To estimate the second gradient of $ \Theta $, we argue  similarly as the above. 
Multiplying both sides of \eqref{3.18} by $ \partial _i \partial _j \Theta | \nabla^2\Theta  |^{ q-2}$,  and taking the sum over $ i,j =1,2$, we find 
\begin{align}
&\frac{1}{q}\partial _t | \nabla ^2 \Theta|^q   + \frac{1}{q}W\cdot \nabla | \nabla^2  \Theta|^q  
\cr
& \qquad= 
 -\sum_{i,j=1}^2 \partial_i W\cdot  \nabla \partial _j\Theta   \partial_i  \partial _j\Theta | \nabla^2\Theta  |^{ q-2}-  
 \sum_{i,j=1}^2\partial_j W \cdot \partial _i\nabla  \Theta \partial_i \partial_j  \Theta | \nabla^2\Theta  |^{ q-2}
\cr
& \qquad\qquad \qquad  -\sum_{i,j=1}^2 \partial_i \partial_j  W\cdot \partial_i \partial_j  \Theta \cdot \nabla \Theta | \nabla^2\Theta  |^{ q-2}
\cr
& \qquad\le 2| \nabla W| | \nabla^2\Theta  |^{ q} +| \nabla^2 W| | \nabla \theta | | \nabla^2\Theta  |^{ q-1}.  
\label{3.70}
\end{align}
Integrating both sides of \eqref{3.70} over $ B(r_0)$, and applying integration by parts, 
using H\"older's inequality, and observing \eqref{3.5},  we infer  
\begin{align}\label{3.70a}\frac{1}{q} \partial _t \| \nabla^2 \Theta (t) \|_{ L^q(B(r_0))}^q 
   &\le   2\| \nabla  W(t)\|_{ L^\infty(B(r_0))} \| \nabla^2 \Theta (t)\|^{ q}_{ L^q(B(r_0))}
\cr   
  &\qquad   + \| \nabla  \Theta (t)\|_{ L^\infty(B(r_0))}
\| \nabla^2 W (t)\|_{ L^q(B(r_0))}   \| \nabla^2 \Theta (t)\|^{ q-1}_{ L^q(B(r_0))}\cr
&\le  ( 2\| \nabla  W(t)\|_{ L^\infty(B(r_0))} +1)\| \nabla^2 \Theta (t)\|^{ q}_{ L^q(B(r_0))}  \cr
  &\qquad +  c\| \nabla  \Theta (t)\|_{ L^\infty(B(r_0))}^q
\| \nabla^2 W (t)\|_{ L^q(B(r_0))}^q,
\end{align}
where we used Young's inequality in the second inequality. 
From \eqref{3.67}  and \eqref{3.69a} combined with \eqref{defw} we find the term of the last line of \eqref{3.70a}  belongs to $L^\infty(t_0, 0)$. 
Since  the function $ \| \nabla W( \cdot )\|_{ L^{ \infty}(B(r_0))}$  
belongs to $ L^\infty(t_0, 0)$, once more using Gronwall's lemma, as above  we see that $ \nabla ^2 \Theta \in  L^1(t_0, 0; L^q(B(r_0)))$.
Again observing \eqref{3.70a}, we get 
\begin{equation}
\partial _t \| \nabla^2 \Theta(t)  \|_{ L^q(B(r_0))}^q\in L^1(t_0, 0),\quad  \Theta \in L^\infty(t_0, 0; W^{2,\, q}(B(r_0))) . 
\label{3.70}
\end{equation}
 Let $ \phi \in C^{\infty}_{\rm c}(B(1))$ denote a cut off function with $ \phi \equiv 1$ on $ B(r)$, $ 0< r< 1$. 
We apply $ \partial _i \partial _j $, $ i,j \in \{ 1,2\}$ to \eqref{bou1} and multiply the resultant equation by $ \phi $. 
This gives 
\[
\partial_t (\partial _i \partial _j \theta  \phi) + a u \cdot \nabla  (\partial _i \partial _j \theta  \phi) 
= f,
\]  
where 
\[
f= -\partial _i (au)\cdot \nabla \partial_j \theta \phi - \partial _j (au)\cdot \nabla \partial_i \theta \phi  
-  \partial _i\partial _j (au)  \nabla \theta \phi   +  a u \cdot \nabla \phi  \partial _i \partial _j \theta. 
\]
Thanks to \eqref{3.69a}, \eqref{3.69b} and \eqref{3.70} we see that $ f, \partial _i \partial _j \theta  \phi \in L^q(\R^{2}\times (t_0,0))$.   
Furthermore, in view of \eqref{3.69a} and Sobolev's embedding theorem we get $ au\in L^\infty(t_0, 0; W^{1,\, \infty}_{ loc}(B(1)))$.  
Thus, we are in a position to apply Lemma\,\ref{lemC.1} with $ h= \partial _i \partial _j \theta  \phi  $ and $ au$ in place of $ u$. 
Accordingly, $ \partial _i \partial _j \theta  \phi \in C([t_0, 0]; L ^q(\R^{2}))$, which implies  
\begin{equation}
\theta \in C([t_0, 0]; W^{2,\, q}(B(r)))\quad  \forall\,0< r< 1. 
\label{3.71}
\end{equation}   
This  together with $ \theta \in   C([-1, t_0]; W^{ 2,q}(B(1)))$ yields 
$ \theta \in   C([-1, 0]; W^{ 2,q}(B(r)))$ for all $ 0<r<1$. 
  By a similar reasoning we infer from \eqref{3.10} that $ w \in   C([-1, 0]; W^{1,q}(B(r)))$ for all $ 0< r< 1$. This together with 
  \eqref{bou4} and \eqref{3.69b} implies $ u\in C([-1, 0]; W^{2,\, q}(B(r)))$ for all $ 0< r< 1$.    
 Whence, the proof of Theorem\,\ref{thm3.1} is complete.
\hfill \Beweisende

\section{Proof of Theorem\,\ref{thm1.1} and \ref{thm1.2}}
\label{sec:-4}
\setcounter{secnum}{\value{section} \setcounter{equation}{0}
\renewcommand{\theequation}{\mbox{\arabic{secnum}.\arabic{equation}}}}

In the proof of Theorem\,\ref{thm1.1} below we make use of the following 

\begin{lem}
\label{lem4.1}
Let $ B(\xi _0, \rho )$  be a ball in $ \R^{2}$. Furthermore, let $ \phi \in C^{\infty}_{\rm c}(B(\xi _0, \rho))$ denote a cut off function 
such that $ 0 \le \phi \le 1$, and $ | \nabla \phi | \le c \rho^{ -1}$. Then for every $ u \in L^2(B(\xi _0, \rho)) $ with 
$ \curl u \in L^\infty(B(\xi _0, \rho))$ it holds 
\begin{equation}
\| u\phi^{ 5} \|_{ \infty} \le c \| u\phi^{ 4} \|^{ \frac{1}{2}}_{ 2}
\| \curl u \phi^{ 4} \|^{ \frac{1}{2}}_{ \infty} + c \rho ^{ -1} \|  u  \phi ^2\|_{ 2}.
\label{4.18}
\end{equation}

\end{lem} 

{\bf Proof}: 
By virtue of   Lemma\,\ref{lemA.9} with $ n=2$
and Lemma\,\ref{lemA.10} we get 
 \begin{align}
\|  u   \phi  ^5 \|_{ \infty} &\le c \|u \phi ^5 \|_{ 2}^{ \frac{1}{2}}   
| \nabla (u\phi ^5)  |_{ BMO}^{ \frac{1}{2}} 
\cr
&\le c \|  u  \phi ^5 \|_{ 2}^{ \frac{1}{2}}   
\Big(| \curl ( u  \phi ^5)  |_{ BMO} + | \nabla \cdot  ( u  \phi^5)  |_{ BMO}\Big)^{ \frac{1}{2}}
\cr
&\le c \|  u  \phi ^5 \|_{ 2}^{ \frac{1}{2}}   
\Big(| \curl u \phi^5 |_{ BMO} + | u \cdot \nabla \phi^5 |_{ BMO}\Big)^{ \frac{1}{2}}
\cr
&\le c \|  u  \phi ^5\|_{ 2}^{ \frac{1}{2}}   
\| \curl u \phi^5 \|_{ \infty}^{  \frac{1}{2} } +c \rho ^{ - \frac{1}{2}} \|  u \phi ^5 \|_{ 2}^{ \frac{1}{2}}   
\| u  \phi^{ 4}  \|_{ \infty}^{ \frac{1}{2}}.
\label{4.19}
\end{align}
Using Calder\'on Zygmund's inequality,  we estimate  
\begin{align}
\| \nabla u \phi^{ 4} \|_{ 4} &\le c  \| \curl  u \phi^4\|_{ 4} +
 \| \nabla \cdot u \phi^4 \|_{ 4} 
+ c \rho ^{ -1}  \|   u  \phi ^{ 3} \|_{ 4} 
\cr
&\le 
 c \| \curl u \phi^4 \|_{ 4} +
  c \rho ^{ -1}  \|   u  \phi^{ 3} \|_{ 4}
  \cr
&\le  c \rho ^{ \frac{1}{2}}\| \curl u \phi^4 \|_{ \infty} +
  c \rho ^{ -1}  \|   u  \phi^{ 3} \|_{ 4}. 
  \label{4.20}
\end{align}
Noting that $ \|   u  \phi^{ 3} \|_{ 4} \le \|   u  \phi ^{ 2} \|_{ 2}^{ \frac{1}{2}}
 \|   u  \phi ^{ 4} \|_{ \infty}^{ \frac{1}{2}} $, we infer from \eqref{4.20}
 \begin{equation}
 \| \nabla u\phi^{ 4} \|_{ 4} \le 
   c \rho ^{ \frac{1}{2}}\| \curl u\phi^4 \|_{ \infty} +
  c \rho ^{ -1}  \|   u  \phi^{ 3} \|^{ \frac{1}{2}}_{ 2} \|   u  \phi^{4} \|^{ \frac{1}{2}}_{ \infty} . 
 \label{4.21}
 \end{equation} 

On the other hand, by Gagliardo-Nirenberg's inequality we find 
\begin{align*}
\|  u  \phi^{ 4} \|_{ \infty} &\le c\|  u  \phi^{ 4} \|_{ 2}^{ \frac{1}{3}}
\| \nabla  u  \phi^{4} \|_{ 4}^{ \frac{2}{3}} + c \rho ^{ - \frac{2}{3}} 
\|  u  \phi^{ 4} \|_{ 2}^{ \frac{1}{3}} \|   u  \phi^{ 3} \|_{ 4}^{ \frac{2}{3}} 
\\
&\le c\|  u  \phi^{ 4} \|_{ 2}^{ \frac{1}{3}}
\| \nabla  u  \phi^{ 4} \|_{ 4}^{ \frac{2}{3}} + c \rho ^{ - \frac{2}{3}} \|  u  \phi^{2} 
\|_{ 2}^{ \frac{2}{3}}  \|   u  \phi ^{ 4} \|_{ \infty}^{ \frac{1}{3}}.
\end{align*}
Applying Young's inequality, we  obtain 
\begin{equation}
\|  u  \phi^{ 4} \|_{ \infty} \le 
c\|  u  \phi^{ 4} \|_{ 2}^{ \frac{1}{3}}
\| \nabla  u  \phi^{ 4} \|_{ 4}^{ \frac{2}{3}} + c \rho ^{ - 1} \|  u  \phi^{2} \|_{ 2}. 
\label{4.22}
\end{equation}
Combining \eqref{4.21} and \eqref{4.22}, and applying Young's inequality, we deduce that     
\begin{equation}
\|  u  \phi^{ 4} \|_{ \infty} \le 
c \rho ^{ \frac{1}{3}}\|  u  \phi^{ 4} \|_{ 2}^{ \frac{1}{3}}
\| \curl u\phi^{ 4} \|_{ \infty}^{ \frac{2}{3}} + c \rho ^{ - 1} \|  u  \phi^{2} \|_{ 2}. 
\label{4.23}
\end{equation}
We estimate the second term on the right-hand side of \eqref{4.19} by using \eqref{4.23}   and then apply Young's inequality. 
This leads to the estimate \eqref{4.18}.  \hfill \Beweisende  

\vspace{0.5cm}  
{\bf Proof of Theorem\,\ref{thm1.1}}:   
Let $ v\in C([-1,0); W^{2,\, q} (\Bbb R^3))$, $ 3<q<+\infty$, be a solution to the Euler equation \eqref{euler} satisfying 
\eqref{type1i} for some ball $ B(x_{ \ast}, R_0)$.

 In our discussion below let $ \xi_{ \ast}  :=  (\sqrt{x_{ 1, \ast}^2 + x_{ 2, \ast}^2},  x_{ \ast, 3})\in \R^{2}$.   
 In order to apply  Theorem\,\ref{thm3.1},  our first aim will be to  check  that for all $ 0< R< R_0$ the following conditions holds
\begin{equation}
 \intl_{-1}^{0}  \| \widetilde{v}  (t)\|_{ L^\infty(B(\xi  _{ \ast}, R))}  dt < +\infty, 
\label{4.5}
\end{equation}  
where $ \widetilde{v} = (v^r, v^3)$,   and $ B(\xi _{ \ast}, R)= 
\{ y\in \R^{2} \,|\, | y- \xi _{ \ast}| <R\}$.

 \vspace{0.2cm}
{\it  Proof of \eqref{4.5}}: 
Observing \eqref{vort-th1},  we see that 
 \begin{equation}
 \partial _t \omega^\theta  + \widetilde{v}\cdot \widetilde{\nabla} \omega^\theta  = \frac{v^r}{r} \omega ^\theta  + 2\frac{v^\theta \partial _3 v^\theta }{r} 
 \quad \text{ in}\quad  \R^{2}_+ \times (-1,0),    
 \label{4.6}
 \end{equation} 
 where 
$ \widetilde{\nabla } = (\partial _r, \partial _3).$

Let $ 0<R<R_0$ be  arbitrarily chosen, but fixed. Set $ R_1 = \frac{R_0+R}{2}$. 
  Let $ \phi \in C^{\infty}_{\rm c}(B(\xi _{ \ast}, R_0))$ such that 
 $ 0 \le \phi \le 1$, $ \phi \equiv 1$ on $ B(\xi _{ \ast}, R_1)$ and $ | \nabla \phi | \le c (R_0-R)^{ -1}$. 
We multiply both sides of \eqref{4.6}  by $ \phi ^6 {\rm sign}(\omega^\theta ) | \omega^\theta |^{ - \frac{1}{2}}$,  and set $ V= \phi ^6 | \omega^\theta  |^{ \frac{1}{2}}$. This gives
 \begin{equation}
 \partial _t V+ \widetilde{v}\cdot \widetilde{\nabla} V = 6\nabla \phi \cdot \widetilde{v}  \phi^5 | \omega ^\theta |^{ \frac{1}{2}}+ \frac{v^r}{r} V + 
 2\frac{ \phi ^{ 12}v^\theta \partial _3 v^\theta }{r{\rm sign}(\omega^\theta  ) V} 
 \quad \text{ in}\quad  \R^{2}_+ \times (-1,0).    
 \label{4.7}
 \end{equation} 
Let  $ (\rho  , z, s)\in B(\xi _{ \ast}, R_1)\times  (-1,0)$ be  arbitrarily chosen. By $ X(\tau ): (-1, s) \rightarrow \R^{2}_+$ we denote the trajectory such that 
\[
\dot{X}(\tau ) = \widetilde{v} (X(\tau ), \tau ),\quad  X(s) = \xi _0 = (\rho , z).  
\] 

We claim that there exists a constant $ c>0$ independent of $ (\rho , z, s)$ such that  
\begin{align}
& | V(\rho , z, s)| \le \max\{ (-s)^{ -  \frac{1}{2} }, \| \omega ^\theta (-1)\|_{\infty}\} 
\cr
&\qquad \qquad \quad + c 
 \intl_{t_0}^{s} \Big( \| \omega ^\theta (\tau )\|_{\infty } +\|\partial _3 v^\theta (\tau ) \|_{ \infty}\Big)d\tau 
+  c. 
\label{claim}
\end{align}

In case $ | V(\rho , z, s)| \le  (-s)^{ - \frac{1}{2}}$, \eqref{claim} is obvious. Thus, without the  loss of generality we may assume that  
 $ V(\rho , z, s) > (-s)^{ - \frac{1}{2}}$.  There are two cases, either $ V(X(\tau ), \tau ) > (-\tau )^{ - \frac{1}{2}}$ 
for all $ \tau \in (-1, s)$, then we set $ t_0=-1$, or there exists $ -1 \le  t_0< s$ such that 
\[
V(X(t_0 ), t_0) = (-t_0)^{ - \frac{1}{2}},\quad  V(X(\tau ), \tau ) \ge  (-\tau )^{ - \frac{1}{2}}\quad  \forall\,\tau \in (t_0, s). 
\]
The sign of $ \omega^\theta  $ does not change in $ [t_0, s]$(since $|\omega^\theta|$ does not touch zero), we may assume that $ {\rm sign}(\omega^\theta  )=1$ in $ [t_0,s]$. 
Using the chain rule, we derive from \eqref{4.7} 
\begin{align}
&\frac{d}{d\tau } V(X(\tau ), \tau ) = F(\tau )+ G(\tau ),
\label{4.9}
\end{align}
where we set
\begin{align*}
F(\tau ) &=6 r(X(\tau ))^{-1} \nabla \phi(X(\tau )) \cdot  \widetilde{v} (X(\tau ), \tau ) \phi(X(\tau ))^{ 5}\omega ^\theta (X(\tau ), \tau )^{ \frac{1}{2}} 
\\
&\qquad \qquad + r(X(\tau ))^{ -1} V(X(\tau ), \tau ) v^r (X(\tau), \tau ), 
\\
G(\tau ) &= 2\frac{ \phi(X(\tau ))^{ 12}v^\theta (X(\tau ), \tau )\partial _3 v^\theta(X(\tau ), \tau ) }{r(X(\tau ))V(X(\tau ), \tau )},
\quad \tau \in (t_0, s). 
\end{align*}

Integrating both sides  of \eqref{4.9} over $ (t_0, s)$, we are led to 
\begin{align}
V(R, z, s) & \le  \max\{ (-t_0)^{ - \frac{1}{2}}, V(X(-1))\} +  \intl_{t_0}^{s} ( F(\tau )+ G(\tau ) )d\tau
\cr
&\le \max\{ (-s)^{ - \frac{1}{2}},  \| \omega ^\theta(-1)\|^{ \frac{1}{2}}_{ L^\infty(\R^{2}_+)}  \} +  \intl_{t_0}^{s}  (F(\tau )+ G(\tau ) )d\tau. 
\label{4.10}
\end{align}

Firstly,  it is readily seen that 
\[
F(\tau ) \le c \| \widetilde{v} (\tau ) \phi ^5\|_{\infty} \| \omega ^\theta(\tau ) \|_{ \infty}^{ \frac{1}{2}}. 
\]
Thanks to \eqref{4.18} with $ u= \widetilde{v} (t)$, using Young's inequality,  and recalling  that 
$ v\in L^\infty(-1, 0; L^2(\R^{3}))$,  we get 
\begin{align}
F(\tau ) &\le c \| \widetilde{v}(\tau ) \phi^{ 4} \|^{ \frac{1}{2}}_{ 2}
\| \omega ^\theta(\tau ) \|_{ \infty} + c \| \widetilde{v}(\tau ) \phi ^2\|_{ 2} \| \omega ^\theta(\tau ) \|^{ \frac{1}{2}}_{ \infty}
\cr
&\le c \| r^{ \frac{1}{2}}\widetilde{v}(\tau )  \|^{ \frac{1}{2}}_{ 2}
\| \omega ^\theta(\tau ) \|_{ \infty} + c \|r^{ \frac{1}{2}} \widetilde{v}(\tau ) \|_{ 2} \| \omega ^\theta(\tau ) \|^{ \frac{1}{2}}_{ \infty}
\cr
&\le c  \| \omega ^\theta(\tau ) \|_{ \infty} + c,
\label{4.11}
\end{align}
where $ c>0$, depending only on $ (R_0-R)^{ -1}$. 

 Secondly, 
 \begin{align}
 | G(\tau )| \le c (-\tau )^{ \frac{1}{2}} \| r v^\theta  \|_{ L^\infty(B(\xi _{ \ast}, R_0))}
  \| \partial _3 v^\theta(\tau ) \|_{L^\infty(B(\xi _{ \ast}, R_0))}
 \le c   \| \partial _3 v^\theta(\tau ) \|_{L^\infty(B(\xi _{ \ast}, R_0))}
 \label{4.12}
  \end{align}
where the constant $ c$ depends only on $ (R_0-R)^{ -1}$. 
 
 Inserting \eqref{4.11} and \eqref{4.12} into \eqref{4.10}, we arrive at 
 \begin{align*}
| \omega ^\theta (\rho , z, s)|^{ \frac{1}{2}} &= V(\rho , z, s) 
\cr
&\le \max\{ (-s)^{ - \frac{1}{2}},  \| \omega ^\theta(-1)\|_{ \infty}  \} 
\\
&\qquad + 
c \intl_{t_0}^{s}  \| \omega ^ \theta (\tau ) \|_{ L^\infty(B(\xi _{ \ast}, R_0))}+ \| \partial _3 v^\theta(\tau ) \|_{L^\infty(B(\xi _{ \ast}, R_0))} d\tau + 
c 
\end{align*}
with a constant $ c>0$ depending on $ (R_0-R)^{ -1}$ but  independent of $ (\rho ,z, s)$. This completes  the proof of \eqref{claim}. 

Accordingly, 
\begin{equation}
\| \omega ^\theta (s)\|_{ L^\infty(B(\xi _{ \ast}, R_1))}^{ \frac{1}{2}} \le c  s^{ - \frac{1}{2}} + c\intl_{t_0}^{s}  \| \omega  (\tau ) \|_{ L^\infty(B(\xi _{ \ast}, R_0))} d\tau, 
\label{4.14}
\end{equation}
with a constant $ c>0$ depending on $(R_0-R)^{ -1}$, but not on $ z$ and  $ s$.  

On the other hand, for given $ (\rho , z, s)\in B(\xi _{ \ast}, R)$, we may choose a cut-off function 
$ \phi \in C^{\infty}_{\rm c}(B((\rho ,z), \sigma ))$, with $ \sigma = \frac{R_0-R}{4}$ such that $ \phi (\rho , z) = 1$ and 
$ | \nabla \phi | \le c\sigma  ^{ -1}$. Then we  apply \eqref{4.18} with $ u= \widetilde{v}(s) $, which shows that 
\begin{align}
| \widetilde{v}(\rho , z, s)| &\le c\| \widetilde{v} (s) \phi ^4\|_{ 2}^{ \frac{1}{2}} 
\| \omega ^\theta (s)\|_{L^\infty(B(\xi _{ \ast}, R_1)}^{ \frac{1}{2}}
+ c \sigma ^{ -1} \| \widetilde{v} (s) \phi ^4\|_{ 2}
\cr
&\le c x_{ 1, \ast}^{ - \frac{1}{4}}\| r^{ \frac{1}{2}}\widetilde{v} (s) \|_{ 2}^{ \frac{1}{2}} \| \omega ^\theta (s)\|_{ L^\infty(B(\xi _{ \ast}, R_1)}^{ \frac{1}{2}}
+ c \sigma ^{ -1} x_{ 1, \ast} ^{ - \frac{1}{2}} \| r^{ \frac{1}{2}}\widetilde{v} (s) \|_{ 2}
\cr
&\le  c\| \omega ^\theta (s)\|_{L^\infty(B(\xi _{ \ast}, R_1)}^{ \frac{1}{2}}
+ c,  
 \label{4.14aa}
\end{align}
Thus, combining  \eqref{4.14aa} and  \eqref{4.14}, we get 
\begin{equation}
\label{4.14a}
\| \widetilde{v}(s) \|_{L^\infty(B(\xi _{ \ast}, R))} \le c (1+\| \omega ^\theta (s)\|_{ L^\infty(B(\xi _{ \ast}, R_1))}^{ \frac{1}{2}})
\le c  s^{ - \frac{1}{2}} + c\intl_{t_0}^{s}  \| \omega  (\tau ) \|_{ L^\infty(B(x _{ \ast}, R_0))} d\tau.
\end{equation}

Since  the right-hand side of \eqref{4.14a} belongs to $ L^1(-1,0)$, we have \eqref{4.5}.

\vspace{0.5cm}  
Let  $0<R< \rho_{ \ast}$,   where $\rho_{ \ast}=\sqrt{x_{1, \ast}^2+x_{2, \ast}^2}$ is given in Theorem\,\ref{thm1.1}.  Setting  
\[
(r, x_3):=\Big(\rho_{ \ast} + R y_1, x_{ 3, \ast} + Ry_2\Big),\quad  (y_1, y_2)\in B(1), 
\]
we define 
\begin{align*}
u (y_1, y_2, t) &= r (v^r(r, x_3, t), v^3(r, x_3, t)) 
\\
\theta (y_1, y_2, t) &= (r v^\theta(r, x_3, t))^2,  
\\
w (y_1, y_2, t) &= \frac{\omega ^\theta(r, x_3, t) }{r},
\end{align*} 
we see that $ (u , \theta  , w  )$ solves the system \eqref{bou1},  \eqref{bou2}, \eqref{bou3}   
in $ B(1)\times (-1,0)$ with 
\[
a(y) = \frac{1}{R(\rho_{ \ast} + R y_1)},\quad b_1(y)=0, \quad b_2(y)   = \frac{1}{R(\rho_{ \ast} + R y_1)^4},
\]
while the equation    $ \omega ^\theta = \partial _3 v ^r - \partial_r v^3$  turns into  
\begin{align*}
\curl u  &:=\partial_1 u_2 -\partial _2 u_1=d(y) w  +e(y)\cdot u
\end{align*} 
in $ B(1)\times (-1,0)$, where 
\[
d(y) =R(\rho _{ \ast}+ R y_1)^2,\quad  e_1(y)=0,\quad   e_2(y) = \frac{R}{\rho_{ \ast} + R y_1}.  
\]

Obviously, $ a,b,d, e\in C^\infty(\overline{B(1)} )$.   Furthermore, by our assumption \eqref{type1i} we get 
 \eqref{type1-bou}$ _1$. 
 Indeed, recalling the relation $ \partial _3 v^\theta = \omega ^r$ and $\partial _r v^\theta = \omega ^3- \frac{v^\theta }{r}$, we see that 
\[
\partial _1 \theta = 2R  rv^{ \theta } \omega ^3(r, x_3),\quad \partial _2 \theta = - 2R   rv^{ \theta } 
\omega ^r(r, x_3).  
\]
Taking  into account of the fact  that   $ r v^\theta $  preserved along the particle trajectories (cf.  \eqref{trans}),  \eqref{type1-bou}$ _1$ 
follows from  \eqref{type1i}. 

Thus $ (u,\theta , w)$ 
solves \eqref{bou3}, \eqref{bou1}, \eqref{bou3}, \eqref{bou4}, and  \eqref{type1-bou}$ _1$ holds.     
In addition thanks to \eqref{4.5} the condition in  \eqref{type1-bou} also holds. 
In order to apply Theorem\,\ref{thm3.1} it only remains to verify that $ u\in C_w([-1, 0]; L^2(B(1)))$.  Indeed, 
recalling that $ v \in L^\infty(-1, 0; L^2(\R^{3}))$, and noting that the following identity 
\[
 \intl_{-1}^{t} \intl_{\R^{3}}  v(s) \cdot \partial _t \varphi (s) + v(s) \otimes v(s) :\nabla \varphi (s)dx ds = \intl_{\R^{3}} v(t) \cdot \varphi (t) dx
\]
holds true for all $ t\in (-1, 0)$ and for all $ \varphi \in C^{\infty}_{\rm c}(\R^{3}\times (-1,0 ])$ with $ \nabla \cdot \varphi =0$.  Therefore there exists a unique $ v(0)\in L^2(\R^{3} )$ such that 
$ v(t) \rightarrow v(0)$ weakly in $ L^2(\R^{3} )$ as $ t \nearrow 0$. By this definition of $ v$ we get   
 $ v\in C_w([-1, 0]; L^2(\R^{3}))$. Thus, by virtue of the definition of $ u$ we get $ u\in C_w([-1,0]; L^2(B(1)))$.
Accordingly, we are in a position to apply Theorem\,\ref{thm3.1}. This completes the proof of Theorem\,\ref{thm1.1}.  \hfill \Beweisende  

\vspace{0.3cm}
{\bf Proof of Theorem\,\ref{thm1.2}}:  Since \eqref{type2} implies \eqref{type1i},  the assertion of  Theorem\,\ref{thm1.2} 
is an immediate consequence of  Theorem\,\ref{thm1.1}.  \hfill \Beweisende

\section{Proof of Theorem\,\ref{thm1.3}}
\label{sec:-5}
\setcounter{secnum}{\value{section} \setcounter{equation}{0}
\renewcommand{\theequation}{\mbox{\arabic{secnum}.\arabic{equation}}}}

Given $R>0$,  we denote below  $H_R =\{ x\in \Bbb R^3\, |\, x_1^2 +x_2 ^2 >R^2 \}$.
Thanks to  Theorem\,\ref{thm1.1}  the statement of Theorem\,\ref{thm1.3} 
will be an immediate consequence of the following.

\begin{lem}
\label{thm2.1} 
Let $ v\in C([-1,0); W^{1,\, \infty}(\R^{3})) \cap  L^\infty(-1,0; L^{2}(\R^{3})),$ 
$ 2<q<+\infty, $  be an axisymmetric solution to  \eqref{euler} in $ \R^{3}\times (-1,0)$. If the for $ 0<R<+\infty$ the following two conditions are  fulfilled 
\begin{equation}
 \intl_{-1}^{0} (-t)   \| \nabla v^\theta   (t) \|_{ L^\infty(H_{ \frac{R}{4}})}   dt <+\infty, \quad 
 \intl_{-1}^{0}  \| v^r  (t)\|_{ L^\infty(H_{ \frac{R}{4}})}   dt   < +\infty. 
\label{type1ii}
\end{equation}
Then
\begin{align}
 \intl_{-1}^{0}  \| \omega^\theta  (t)\|_{ L^\infty(H_R)}  dt < +\infty,
\label{bkm-vor}
\\
\limsup_{t \to 0^-} (-t) \| \omega^\theta  (t)\|_{ L^\infty(H_R)} =0.
\label{type1-vor}
\end{align}
\end{lem}

{\bf Proof}:  Let $ 0<R< +\infty$ be fixed.  
According to \eqref{type1ii}$ _2$ there exists $ t_0\in (-1, 0)$, such that 
\begin{equation}
\intl_{t_0}^{0}  \| v^r(s) \|_{L^\infty(H_{ \frac{R}{4}}) } ds\le \frac{R}{16}. 
\label{def_t0}
\end{equation}
 We now set  $ \varrho (r, t)= r + 8R^{ -1}(r-R)  \intl_{t}^{0}  \| v^r(s) \|_{L^\infty(H_{ \frac{R}{4}}) } ds$.  Owing to \eqref{def_t0} 
 we see that for all $t\in (t_0, 0)$
 \begin{equation}
\varrho (r, t) \ge r + \frac{r-R}{2} = \frac{3r -R}{2} \ge \frac{R}{4}\quad \forall\,\frac{R}{2} \le r <+\infty. 
 \label{2.2}
 \end{equation}
 
 Clearly,
 \begin{align*}
\partial _t \varrho (r, t) &=  - 8R^{ -1}\| v^r(s) \|_{L^\infty(H_{ \frac{R}{4}}) } (r-R),\quad  
  \\
 \partial _r \varrho (r, t) &= 1+ 8R^{ -1}\intl_{t}^{0}  \| v^r(s) \|_{L^\infty(H_{ \frac{R}{4}}) } ds. 
 \end{align*}
We define 
 \begin{align*}
V (r,x_3,t) &= v(\varrho (r, t), x_3, t),
\\
\Theta(r,z,t)  &=  v^\theta  (\varrho (r, t), x_3, t),
\\
\Omega(r,x_3,t) &=\varrho (r, t)^{ -1} \omega^{ \theta }  (\varrho (r, t), z, t), \quad  (r,z,t)\in H_{ \frac{R}{2}}\times (t_0, 0). 
\end{align*}
In addition, define  
\[
\begin{cases}
W^r(r,x_3 ,t) =  \frac{-\partial _t \varrho (r, t) + V^r(r,x_3,t)}{\partial _r \varrho(r, t) },
\\[0.3cm]
W^3(r, x_3 ,t) =   V^3(r,x_3,t),\quad  (r,x_3,t)\in H_{ \frac{R}{2}}\times (t_0, 0). 
\end{cases}
\]

Note that, \eqref{2.2} implies  
\begin{equation}
\varrho (r,t) \ge  \frac{r}{2} \quad  \forall\, \frac{R}{2} < r< +\infty, t_0 < t \le  0. 
\label{2.2a}
\end{equation}

Therefore the functions $ V, \Theta $ and $ \Omega $ are well defined on $ H_{ \frac{R}{2}}\times (t_0, 0).$ 

\hspace{0.5cm}
To proceed, we verify  
\begin{equation}
 W^r(r,x_3,t) \le  -\frac23\| v^r(t) \|_{ L^\infty(H_{ \frac{R}{4}})}\quad  \forall\, (r,x_3,t)\in (\overline{H} _{ \frac{R}{2}}  \setminus  H _{ \frac{3R}{4}}) \times (t_0, 0). 
\label{2.5}
\end{equation}
In fact, according to \eqref{2.2},   we estimate   
\begin{align*}
\partial _r \varrho (r, t)W^r(r,x_3,t)&=
- \partial _t\varrho  (r,t) + V^r(r,x_3, t)
\\
&=  8 R^{ -1}(r-R) \| v^r(t) \|_{ L^\infty(H_{ \frac{R}{4}})} + V^r(r,x_3, t)
\\
& \le  -2 \| v^r(t) \|_{L^\infty( H_{ \frac{R}{4}})} + \| V^r(t)\|_{ L^\infty(H_{ \frac{R}{2}})}
\\
&\le -2 \| v^r(t) \|_{ L^\infty(H_{ \frac{R}{4}})} + \| v^r(t)\|_{ L^\infty(H_{ \frac{R}{4}})}. 
\end{align*}
Since $1\leq  \partial _r \varrho (r, t) \leq \frac32$ by \eqref{def_t0} it follows \eqref{2.5}.  

\vspace{0.3cm}
By the chain rule we see that \eqref{vort-th} turns into the following equations  in  $ H_{ \frac{R}{2}}\times (t _0,0)$
\begin{equation}
\partial _t \Omega + W\cdot \nabla \Omega = 2\frac{\Theta\partial _3 \Theta }{\varrho(r, t)^2}.  
\label{vort-trans}
\end{equation}
 
For $ (r, x_3, s)\in H_{ \frac{R}{2}}\times (t_0, t)$  by  $X= (X^r, X^3)= X(r, x_3,s; \cdot ): [t_0, s] \rightarrow  \R^{2}$ we denote the particle trajectory  such that 
\[
\dot{X}(\tau ) = W(X(\tau ), \tau ),\quad  X(s) = (r, x_3). 
\]

Let $ (r,x_3, s) \in H_{ \frac{R}{2}}\times (t_0,0)$. 

\vspace{0.3cm}
\hspace{0.5cm}
 We claim that $ X(\tau ) \in H_{ \frac{R}{2}}$ for all $ \tau \in [t_0, s]$. 
Otherwise,  there exists $ \tau _0\in [t_0, s]$ such that $ X(\tau _0)\in \partial H_{ \frac{R}{2}}=\{ r= \frac{R}{2}\}$ and $ X(\tau )\in H_{ \frac{R}{2}}$ for all $ \tau \in (\tau_0, s]$. 
This gives  
\[
\dot{X^r}(\tau _0)  = W^r(X(\tau _0), \tau _0) \ge 0,
\]
which contradicts \eqref{2.5},  since $ \| v^r(t) \|_{ L^\infty(H_{ \frac{R}{4}})}$ must be strictly positive.  Thus, the claim is proved.  On the other hand, by the chain rule,  \eqref{vort-trans} gives 
\[
 \frac{d}{d\tau }  \Omega (X(\tau ), \tau ) = 2\frac{\Theta (X(\tau ), \tau )\partial _3 \Theta(X(\tau ), \tau )}{\varrho(X^r(\tau ), \tau )^2}. 
\]   
Recalling that $ X(s)=(r, x_3)$, integration over $ (s_0, s), t_0 \le s_0 \le t$, yields 
\[
\Omega (r, x_3, s)= \Omega (X(s_0), s_0) +  2\intl_{s_0}^{s}  \frac{\Theta(X(\tau ), \tau )\partial _3 \Theta(X(\tau ), \tau )}{\varrho (X^r(\tau ), \tau )^2} d\tau.
\]
Accordingly, using \eqref{2.2a},   and 
noting that $ | \varrho (r, t)\Theta (r, x_3, t)|  =  |\varrho (r, t)v^\theta  (\varrho (r, t), x_3, t)|\le  \| r v^\theta \|_\infty \le c <+\infty$,  we get 
\begin{align}
\| \Omega (s)\|_{ L^\infty(H_{ \frac{R}{2}})} &\le \| \Omega (s_0)\|_{ L^\infty(H_{ \frac{R}{2}})} 
+ c R^{ -3} \intl_{s_0}^{s} \| \partial_3   \Theta(\tau )\|_{ L^\infty(H_{\frac{R}{2}})} d\tau . 
\cr
& \le \| \Omega (s_0)\|_{ L^\infty(H_{ \frac{R}{2}})} 
+ c R^{ -3} \intl_{s_0}^{s} \|   \partial _3  v^\theta (\tau  )\|_{ L^\infty(H_{\frac{R}{4}})} d\tau. 
\label{2.9}
\end{align}
In \eqref{2.9} we take   $ s_0=t_0$,  and integrate both sides  over $ (t_0, t)$. This leads to 
\begin{align*}
\hspace{-2cm} \intl_{t_0}^{t}  \| \Omega (s)\|_{ L^\infty(H_{\frac{R}{2}})}  ds
 &\le  (t-t_0 ) \| \Omega (s_0)\|_{ L^\infty(H_{\frac{R}{2}})} 
- c R^{ -3}\intl_{t_0}^{t} (-s)'\intl_{s_0}^{s} \| \partial_3   v^\theta (\tau )\|_{ L^\infty(H_{\frac{R}{4}})} d\tau ds
\\
& \le  (t-t_0) \| \Omega (t_0)\|_{ L^\infty(H_{\frac{R}{2}})} 
+c R^{ -3}\intl_{t_0}^{t} (-s) \|  \partial_3   v^\theta (s)\|_{ L^\infty(H_{ \frac{R}{4}})}  ds. 
\end{align*}
Observing \eqref{type1ii}$ _1$, this proves that \eqref{bkm-vor}. To verify \eqref{type1-vor},  we first note that  \eqref{2.9} multiplied by $ (-s)$ implies   
\begin{equation}
(-s)\| \omega (s)\|_{L^\infty(H_{R})} \le (-s)\| \omega (s_0)\|_{ L^\infty(H_{ \frac{R}{4}})} 
+cR^{ -3}\intl_{s_0}^{0} (-\tau )\|  \partial_3  v^\theta(\tau )\|_{ L^\infty(H_{ \frac{R}{4}})} d\tau. 
\label{2.10}
\end{equation}
Applying $ \limsup$ as $ s \rightarrow 0^-$ to both sides of \eqref{2.10}, we are led to 
\begin{align*}
&\limsup_{ s\to 0^-}(-s)\| \omega (s)\|_{ L^\infty(H_{ R})} 
\\
&= \limsup_{ s\to 0^-}(-s)\| \omega (t_0)\|_{ L^\infty(H_{ \frac{R}{4}})} +cR^{ -3}\intl_{s_0}^{0} (-\tau )
\|  \partial_3 v^\theta(\tau )\|_{ L^\infty(H_{ \frac{R}{4}})} d\tau 
\\
& \le cR^{ -3}\intl_{s_0}^{0} (-\tau )
\| \partial_3 v^\theta(\tau )\|_{ L^\infty(H_{ \frac{R}{4}})} d\tau.
\end{align*}
On the other hand, the integral on the  right-hand side of this inequality tends to $ 0$ as $ s_0 \rightarrow 0$ we  conclude  \eqref{type1-vor}.  
 \hfill \Beweisende

\appendix
\section{Gagliardo-Nirenberg inequalities with cut-off}
\label{sec:-A}
\setcounter{secnum}{\value{section} \setcounter{equation}{0}
\renewcommand{\theequation}{\mbox{A.\arabic{equation}}}}

Below we denote the H\"{o}lder conjugate of $q$ by $q'=\frac{q}{q-1}$.
Following similar argument to the proof of \cite[Lemma\,2.3]{cha5}, we get 
\begin{lem}
\label{lem3.4}
 Let  $ \psi \in C^{\infty}_{\rm c}(B(r))$, $ 0< r< +\infty$, such that $ 0 \le \psi \le 1$ in $ B(r)$. For all $ u\in W^{1,\, q}(B(r))$, $ 2< q < +\infty$ with $ \nabla \cdot u=0$ a.e. in $ B(r)$ and 
 for  all $ m \ge \frac{2}{q'}$ it holds 
 \begin{align}
 \| \nabla  u \psi ^{ m}\|_{ q} &\le c  \| \curl  u \psi^m \|_{ q} 
 +c \| \nabla \psi \|_{ \infty}^{ \frac{2(q-1)}{q}} \| u \psi ^{ m- a}\|_{ 2},
 \label{3.30}
 \\
 \| u \psi ^{ m-k}\|_{ q} &\le c \| u \psi ^{ m- \frac{2k}{q'}}\|_{ 2}^{ \frac{q'}{2}} \| \curl  u \psi^m \|_{ q} ^{ 1- \frac{q'}{2}}
 + c\| \nabla \psi \|_{ \infty}^{ \frac{2}{q'}-1} \| u \psi ^{ m- ka}\|_{ 2}.
 \label{3.31}
 \end{align}
 
\end{lem}

\begin{lem}
\label{lem3.5}
 Let  $ \psi \in C^{\infty}_{\rm c}(B(r))$, $ 0< r< +\infty$, such that $ 0 \le \psi \le 1$ in $ B(r)$. For all $ u\in W^{1,\, q}(B(r))$, $ 2< q < +\infty$ for  all $ m \ge \frac{2q-2}{q}$ it holds 
\begin{equation}
\| u \psi ^{ m}\|_{ \infty} \le c   \| u \psi ^{ m- \frac{q}{2}}\|_{ 2}^{ 1-\frac{q'}{2}} \| \nabla   u \psi^m \|_{ q} ^{ \frac{q'}{2}}
 + c\| \nabla \psi \|_{ \infty} \| u \psi ^{ m-  \frac{q}{2}}\|_{ 2}.
\label{3.32}
\end{equation}
\end{lem}

For further discussion  below we recall the notion of the local BMO space. 
For $ 0<r<+\infty$ we say $ u\in BMO(B(r))$ if 
\[
| u|_{ BMO(B(r))}:=  \sup_{\substack{ x\in B(r)\\0< \rho \le  r} } r^{ -2} 
\intmw_{B(x, \rho )\cap B(r)}  | u(x) - u_{ B(x, \rho )\cap B(r)}| dx   < +\infty,
\]   
Here we have used the following notation for the mean for a given set $ \Omega \subset \R^{2}$ and $ v\in L^1(\Omega )$
\[
v_{ \Omega } = \intmw_{\Omega } v  dx = \frac{1}{m(\Omega )} \intl_{\Omega } v  dx,  
\] 
 where $ m$ stands for the two dimensional Lebesgue measure.

\begin{lem}
\label{lem3.6}
Let    $ \psi \in C^{\infty}_c(B(r))$, $ 0< r< +\infty$, with $ 0 \le  \psi \le 1$. 
 For every  $ u \in W^{1,\, 1}(B(r))$ such that  $ \nabla \cdot u =0$  and $ \curl u\in BMO(B(r))$
it holds
\begin{equation}
\| \nabla  u\psi ^{ 6}\|_{ BMO} \le c \| (\curl u) \psi^{ 5} \|_{ BMO} + 
c  \Big\{r \| \nabla \psi \|^{ 3}_{ \infty}+   \| \nabla \psi \|^{ 2}_{ \infty}\Big\}\| u \psi \|_{ 2}.
\label{3.33}
\end{equation}

\end{lem} 
 
 {\bf Proof}:   Using  Calder\'on-Zygmund's inequality, we find that 
 \begin{align}
\| (\nabla  u) \psi ^{ 6}\|_{ BMO} &\le c\| \curl (\psi^{ 6} u)\|_{ BMO} 
  + c\| u \cdot \nabla \psi^{ 6} \|_{ BMO} 
  \cr
  & \le 
c\| (\curl   u) \psi ^{ 6}\|_{ BMO} + c\| \nabla \psi \|_{ \infty}\| u \psi^{5} \|_{ \infty}. 
\label{3.34}
 \end{align}
 On the other hand, in view of  \eqref{3.32} with  $ q = 4$  and $ m= 5$ we get 
   \begin{equation}
\| u \psi ^{ 5}\|_{ \infty} \le c\| u \psi \|_{ 2}^{ \frac{1}{3 } }   \| \nabla  u \psi ^{ 5}\|^
{ \frac{2}{3}}_{ 4} + 
c\| \nabla \psi \|_{ \infty} \| u \psi\|_{ 2}. 
\label{3.35}
\end{equation}
We estimate the  first term  on the right-hand side of \eqref{3.35}  by \eqref{3.30} with $ q= 4$ and $ m=5$. This together with 
\cite[Lemma\,B.3]{cha5}  gives 
 \begin{align}
\| u \psi ^{ 5}\|_{ \infty} &\le c\| u \psi\|_{ 2}^{ \frac{1}{3} }   \| (\curl  u) \psi ^{ 5}\|^
{ \frac{2}{3}}_{ 4} + c\| \nabla \psi \|_{ \infty}\| u \psi \|_{ 2}. 
 \cr
 &\le c r^{ \frac{1}{3}}\| u \psi\|_{ 2}^{ \frac{1}{3 } }   \|( \curl  u) \psi ^{ 5}\|^
{ \frac{2}{3}}_{ BMO(B(r))} + c\| \nabla \psi \|_{ \infty} \| u \psi \|_{ 2}. 
 \label{3.36}
 \end{align}
Finally, combining \eqref{3.34} and \eqref{3.36}, and applying Young's inequality, we obtain \eqref{3.33}. 

\hfill \Beweisende   

\vspace{0.3cm}
Using the well known John-Nirenberg inequality, we can get the following  
 
\begin{lem}
\label{lemB.3}
Let $ u\in BMO(B(r))$. Then $ u\in \cap_{ 1 \le  q< \infty} L^q(B(r))$, and it holds 
\begin{equation}
\| u\|_{ L^q(B(r))} \le c r^{ \frac{n}{q}} | u|_{BMO(B(r))}+ c r ^{ \frac{n}{q}-n}\| u\|_{ L^1(B(r))} =
c r^{ \frac{n}{q}} \| u\|_{ BMO(B(r))}. 
\label{B.3}
\end{equation}
\end{lem}
For an elementary proof see \cite[Lemma\,B.3]{cha5}.  

\vspace{0.3cm}
\hspace{0.5cm}
Arguing as in \cite{cha5}  we shall show the following.

\begin{lem}
\label{lem3.7} Let $ u \in W^{1,\, 1}(B(r))$ with $\curl u \in BMO(B(r))$. Then for all $ \psi \in C^{\infty}_{\rm c}(B(r))$ 
with $ 0 \le \psi \le 1$ we get 
\begin{align}
&\| (\curl u) \psi ^{5}\|_{ BMO} 
\cr
&\le c \Big\{1+ r^{ 3}\| \nabla \psi \|_{ \infty}^{ 3}\Big\} (| \curl u|_{ BMO(B(r))} + cr^{- 2} \| u\|_{ L^2(B(r))}). 
\label{3.37} 
\end{align}

\end{lem}

{\bf Proof}:
Assume $ r=1$.  Let $ \eta \in C^{\infty}_{\rm c} (B(1))$ such that $ | \nabla \eta | \le c$ and $ \intl_{B(1)} \eta dx \ge c$, where 
$ c >0$ stands for a constant depending only on $ n$. For $ f\in L^1(B(1))$ we define the mean 
\[
\widetilde{f}_{ B(1)} = \frac{1}{ \intl_{B(1)} \eta dx } \intl_{B(1)} f\eta dx.  
\]
First we see that 
\begin{align}
\| \curl u\|_{ L^1(B(1))} &= \| \curl u -\widetilde{ \curl u}_{ B(1)}  \|_{ L^1(B(1))} + 
| \widetilde{ \curl u}_{ B(1)}|
\cr
&\le c| \curl u|_{ BMO(B(1))} + c \| u\|_{ L^1(B(1))}.
\label{3.38}
\end{align}

Using \eqref{3.38}, we  estimate  for $ \rho \ge \frac{1}{2}$ and $ x_0\in \R^{3}$
\begin{align*}
&\intmw_{B(x_0, \rho )} | (\curl u) \psi^5 - (\curl u \psi^5)_{ B(x_0, \rho )}|  dx
\\
&\qquad \le c  \| \curl u\|_{ L^1(B(1))} 
\\
&\qquad \le 
c| \curl u|_{ BMO(B(1))} + c \| u\|_{ L^1(B(1))}.
\end{align*}
In case $ \rho \le \frac{1}{2}$ and $ B(x_0, \rho ) \cap B(1) \neq \emptyset$ 
there exists $ y_0 \in B(1)$ such that $ B(x_0, \rho ) \subset B(y_0, 2\rho )$ and 
\begin{align*}
&\intmw_{B(x_0, \rho )} | \curl u \psi^5 - (\curl u \psi^5)_{ B(x_0, \rho )}|  dx
\\
&\le c\intmw_{B(y_0, 2\rho )} \intmw_{B(y_0, 2\rho )} | \curl u(x) \psi^5(x) - | \curl u(y) \psi^5(y) 
|  dxdy
\\
&\le c| \curl u|_{ BMO(B(1))} + c \intmw_{B(y_0, 2\rho )} \intmw_{B(y_0, 2\rho )} | \curl u(x)  
|  | \psi^5(x) - \psi^5(y)|dxdy.
\end{align*}
By the fundamental theorem of differentiation  and integration we calculate 
\begin{align*}
\psi ^5(y) -\psi ^5(x) &\le 5  \psi ^{ 4}(\xi_1 )\nabla \psi (\xi_1 )\cdot (y-x) 
\\
&= 5\psi ^{ 4}(x)\nabla \psi (\xi _1)\cdot (y-x) +5  (\psi ^{ 4}(\xi_1 )-\psi ^{ 4}(x))\nabla \psi (\xi_1 )\cdot (y-x) 
\\
&= 5\psi ^{ 4}(x)\nabla \psi (\xi _1)\cdot (y-x)
\\
&\qquad + \psi ^{ 2}(\xi _2) \prod_{ i=1}^2 (6-i) \nabla \psi (\xi _i)\cdot (\xi _{ i-1}-x).
\end{align*}
For some $ \xi _i \in [x,y]$, $ i=1, 2$. 
This along with \eqref{3.38}  yields 
\begin{align*}
&\intmw_{B(y_0, 2\rho )} \intmw_{B(y_0, 2\rho )} | \curl u(x)  
|  | \psi^5(x) - \psi^5(y)|dxdy
\\
&\qquad \le c \| \nabla \psi \|_{ \infty}\rho^{ -1}  \intl_{B(y_0, 2\rho )} | \curl u(x)  | \psi^{ 4} (x) dx 
\\
&\qquad \qquad \qquad + c \| \nabla \psi \|_{ \infty}^2 \| \curl u\|_{ L^1(B(1))}
\\
&\qquad \le c \| \nabla \psi \|_{ \infty} \| (\curl u) \psi ^{ 4}\|_{2} 
+ c \| \nabla \psi \|_{ \infty}^2\| \curl u\|_{ L^1(B(1))}. 
\end{align*}
By using H\"older's inequality,  we find that 
\begin{align*}
\| (\curl u) \psi ^{ 4}\|_{ 2}  &\le \|\curl u\|_{ L^{1}(B(1))}^{ \frac{1}{3}}  
\| (\curl u) \psi ^{6}\|_{ 4}^{ \frac{2}{3}}
\\
&\le \|\curl u\|_{ L^{1}(B(1))}^{ \frac{1}{3}}  
\| (\curl u) \psi ^{5}\|_{4}^{ \frac{2}{3}}.
\end{align*}
Applying the embedding $ L^{ 6}(B(1)) \hookrightarrow  BMO(B(1))$ (cf. Lemma\,\ref{lemB.3}), we get 
\[
\| (\curl u) \psi ^{ 4}\|_{2} \le c \|\curl u\|_{ L^{1}(B(1))}^{ \frac{1}{3}}  
\| (\curl u) \psi ^{5}\|_{ BMO}^{ \frac{2}{3}}.
\]

Combining the above inequalities, and applying Young's inequality together with \eqref{3.38}, we arrive at 
\begin{align}
&\| (\curl u) \psi ^{5}\|_{ BMO} 
\cr
& \le c | \curl u|_{ BMO(B(1))}+ c \| \nabla \psi \|_{ \infty}^3 \| \curl u\|_{ L^1(B(1))}
\cr
&\le c (1+ \| \nabla \psi \|_{ \infty}^{ 3}) \Big(| \curl u|_{ BMO(B(1))} + c\| u\|_{ L^1(B(1))}\Big). 
\label{3.39} 
\end{align}
Whence, \eqref{3.37} follows immediately from \eqref{3.39} by  standard scaling argument.  

\hfill \Beweisende

Combining Lemma\,\ref{lem3.7} and Lemma\,\ref{lem3.6}, we get 

\begin{cor}
\label{cor3.8} 
Let    $ \psi \in C^{\infty}_c(B(r))$, $ 0< r< +\infty$, with $ 0 \le  \psi \le 1$.
 For every  $ u \in W^{1,\, 1}(B(r))$ such that  $ \nabla \cdot u =0$  and $\curl  u\in BMO(B(r))$
it holds
\begin{align}
\| \nabla u \psi ^{ 6}\|_{ BMO(B(r))} 
&\le c\Big(1+ r^{ 3}\| \nabla \psi \|_{ \infty}^{ 3}\Big)| \curl u|_{ BMO(B(r))} 
\cr
&  \qquad +c \Big(r^{ -2} + r\| \nabla \psi \|^{ 3}_{ \infty}\Big)\| u\|_{ L^2(B(r))}.
\label{3.40}
\end{align}
\end{cor} 

{\bf Proof}: Combining \eqref{3.33} and \eqref{3.37} along with Young's inequality, we infer 
\begin{align}
&\| \nabla u \psi ^{ 6}\|_{ BMO(B(r))} 
\cr
&\le c \Big(1+ r^{ 3}\| \nabla \psi \|_{ \infty}^{ 3}\Big) (| \curl u|_{ BMO(B(r))} + cr^{ -2} \| u\|_{ L^2(B(r))})
\cr
&  \qquad\qquad  +c  \Big(r \| \nabla \psi \|^{ 3}_{ \infty}+   \| \nabla \psi \|^{ 2}_{ \infty}\Big)\| u \|_{ L^2(B(r))}
\cr 
&  \le c \Big(1+ r^{ 3}\| \nabla \psi \|_{ \infty}^{ 3}\Big)| \curl u|_{ BMO(B(r))} 
   + c \Big(r^{ -2} + r\| \nabla \psi \|^{ 3}\Big)\| u\|_{ L^2(B(r))}.
\label{3.41}
\end{align}
Whence, \eqref{3.40}.  \hfill \Beweisende

\begin{lem}
\label{lem3.9}
Let $ u \in W^{2,\, q}(B(r ))$, $ 2 \le q < +\infty$. Let $ m,k\in \R$ such that  $2 \le   m <+\infty$ and $ 0 < k \le 2m$. 
Then for every $ \psi \in C^{\infty}_c(B(r))$ it holds  
\begin{equation}
\| \nabla u \psi^m \|_{ q} 
 \le c\| u \psi ^{ 2m-k} \|_{ q}^{  \frac{1}{2} } \|\nabla^2  u\psi ^{k}\|_{ q}^{  \frac{1}{2} } +
  c \| \nabla \psi \|_{ \infty}  \| u \psi ^{ m-1}\|_{ q}. 
\label{3.42}
\end{equation}
If in addition, if $ \nabla \cdot u=0$ almost everywhere in $ B(r )$, then for $2  \le  m < +\infty$ and $ m+1 \le k \le 2m$
it holds
\begin{align}
\| \nabla^2 u \psi^k \|_{ q} 
& \le c \|(\nabla \curl  u) \psi ^{k}\|_{ q} + c \| \nabla \psi \|^2_{ \infty} 
 \| u \psi ^{ k-2}\|_{ q}, 
\label{3.43}
\\
\| \nabla u \psi^m \|_{ q} 
 &\le c\| u  \psi^{ 2m-k}\|_{ q}^{  \frac{1}{2} } \|(\nabla \curl   u) \psi ^{k}\|_{ q}^{  \frac{1}{2} } + 
 c \| \nabla \psi \|_{ \infty}  \| u\psi ^{ m-1}\|_{ q}. 
\label{3.44}
\end{align}

\end{lem}

{\bf Proof}:  
Applying integration by parts, and using H\"older's inequality, we get 
\begin{align*}
\|\nabla  u \psi^m \|^{ q}_{q}  
& =-\intl_{B(\rho )} u \nabla u\cdot  \nabla  | \nabla u|^{ q-2} \psi ^{ mq} dx- 
\intl_{B(\rho )} u \nabla u  | \nabla u|^{ q-2} \cdot \nabla \psi ^{ mq} dx
\\
& \le c \| u \psi ^{ 2m-k}\|_{q}  \| \nabla ^2 u \psi ^{ k}\|_{q}  \|\nabla  u \psi^m \|^{ q-2}_{q}  
\\
&\qquad + c \| \nabla \psi \|_{ \infty} \| u \psi^{ m-1}\|_{L^q(B(\rho ))}   \|\nabla  u \psi^{ m} \|^{ q-1}_{q},  
\end{align*}  
and Young's inequality gives \eqref{3.42}.   

\hspace{0.5cm}
Suppose  $ \nabla \cdot u =0$ almost everywhere in $ B(r )$.  We first  apply \eqref{3.30} with 
$ \nabla u$ in place of $ u$ and $ m=k$, and then use  \eqref{3.42}  with $ m= k-1$. This gives 
\begin{align*}
\| \nabla ^2 u \psi^k \|_q &\le c \| (\nabla  \curl u )\psi^k \|_q + 
c\| \nabla \psi \|_{ \infty}\| \nabla u\psi^{ k-1} \|_{ q} 
\\
&\le  c \| (\nabla  \curl u) \psi^k\|_q + c\| \nabla \psi \|_{ \infty}\| u \psi^{ k-2} \|^{ 1/2}_{q} 
\| \nabla^2 u\psi ^{ k}\|^{ 1/2}_{ q}  
\\
&\qquad \qquad + \| \nabla \psi \|_{ \infty}^2\| u \psi^{ k-2}\|_{ q}. 
\end{align*}
Applying Young's inequality,  we obtain \eqref{3.43}.  
The estimate \eqref{3.44} is now an immediate consequence of \eqref{3.42} 
and \eqref{3.43}.   
\hfill \Beweisende  

\vspace{0.3cm}
Combining Lemma\,\ref{lem3.4} and Lemma\,\ref{lem3.9}, we get the following  

\begin{cor}
\label{cor3.10} 
For all $ u \in W^{2,\, q}(B(r))$, for all $ \psi \in C^{\infty}_{\rm c}(B(r))$ with $0 \le \psi  \le 1$ and for all $ k >5$ we get 
\begin{equation}
\| \nabla ^2 u \psi ^k\|_{ q} \le c \| (\nabla \curl u) \psi ^k\|_{ q} +  c\| \nabla \psi \|_{ \infty}^{ 1+ \frac{2}{q'}} 
\| u \psi ^{ k -4}\|_{2}.
\label{3.45}
\end{equation} 
\end{cor}

{\bf Proof}: Let $ k >q$. The estimate \eqref{3.31}  with $ m=k$ and $ k=2$ reads 
\[
\| u \psi ^{ k-2}\|_{ q} \le c \| u \psi ^{ k-  \frac{4}{q'}}\|_{ 2}^{ \frac{q'}{2}} \| \nabla\curl  u \psi^{ k} \|_{ q} ^
{ 1- \frac{q'}{2}}
 + c\| \nabla \psi \|_{ \infty}^{ \frac{2}{q'}-1} \| u \psi ^{ k- \frac{4}{q'}}\|_{ 2}.
\]
Combining this inequality with \eqref{3.43}, and applying  Young's inequality, we obtain \eqref{3.45}.  \hfill \Beweisende   

\vspace{0.3cm}
Next, we shall   establish an Gagliardo-Nirenberg inequality involving the $ BMO$ norm of the gradient, which improves the case with corresponding $L^\infty$ norm.

\begin{lem}
\label{lemA.9}
For $ u\in L^2(\R^{n}) \cap W^{1,\, 1}_{ loc}(\R^{n})$ such that $ \nabla u \in BMO$ it holds 
\begin{equation}
\| u\|_{ \infty} \le c \| u\|_{2}^{ \frac{2}{n+2}} | \nabla u|_{ BMO}^{ \frac{n}{n+2}}. 
\label{A.17}
\end{equation}
\end{lem}

{\bf Proof}:  First let us show that for every  $ u\in L^2(\R^{n}) \cap W^{1,\, 1}_{ loc}(\R^{n})$ such that $ \nabla u \in BMO$ 
it holds 
\begin{equation}
\| u\|_{ \infty} \le c (\| u\|_{ 2} + | \nabla u|_{ BMO})
\label{A.18}
\end{equation} 
for some constant $ c>0$ depending only on $ n$.  To prove \eqref{A.18} let  $ x_0\in \R^{n}$ be fixed.  
We take  a cut off function 
 $ \zeta \in C^{\infty}_{\rm c} (Q(x_0,1))$  such that $ 0 \le \zeta \le 1$, $ | \nabla \zeta | \le 3$,  and  $ \intl_{Q(x_0,1)} \eta dx \ge 2^{ -n}$, and  define the generalized mean 
\[
\widetilde{\nabla u}_{ B(x_0,1)} = \frac{1}{ \intl_{Q(x_0,1)} \zeta dx} \intl_{Q(x_0, 1)} \nabla u \zeta dx.   
\]
By virtue of Sobolev embedding theorem, and Jensen's inequality we find 
\begin{align*}
&\| u - \widetilde{\nabla u}_{ Q(x_0,1)}(x-x_0)\|_{ L^\infty(Q(x_0,1))} 
\\
&\quad\le    c\|  u - \widetilde{\nabla u}_{ Q(x_0,1)}(x-x_0)\|_{ L^2(Q(x_0,1))}  +  c\|  \nabla u - \widetilde{\nabla u}_{ Q(x_0,1)}\|_{ L^{ 2n}(Q(x_0,1)) } 
\\
&\quad  \le c\|  u \|_{ L^2(Q(x_0,1))}+ c |  \widetilde{\nabla u}_{ Q(x_0,1)}|  +  
c\|  \nabla u - (\nabla u)_{ Q(x_0,1)}\|_{ L^{ 2n}(Q(x_0,1)) }.  
\end{align*}
Here, $ Q(x_0, 1)$ stands for the usual cube. Using integration by parts, we see that 
\[
|  \widetilde{\nabla u}_{ Q(x_0,1)}| \le 2^n \bigg|\intl_{Q(x_0, 1)}  u \nabla \zeta dx\bigg| \le c \| u\|_{ L^2(Q(x_0,1))}. 
\]
Furthermore, employing John-Nirenberg's inequality\cite{joh},  we get 
\[
\|  \nabla u - \widetilde{\nabla u}_{ Q(x_0,1)}\|_{ L^{ 2n}(Q(x_0,1)) } \le c| \nabla u|_{ BMO(Q(x_0, 1))} \le 
c| \nabla u|_{ BMO}. 
\]
Combining the last two inequalities,  we obtain the desired estimate  \eqref{A.18}.  

\hspace{0.5cm}
Next, given $ \lambda >0$, we define $ u_\lambda (x) = u(\lambda x), x\in \R^{n}$. Applying the  chain rule together with the 
transformation formula of the Lebesgue integral we find   for any ball $ B(x, r) \subset \R^{n}$
\begin{align*}
&\frac{1}{| B(x, r)|}\intl_{B(x, r)}  | \nabla u_\lambda - (\nabla u_\lambda)_{ B(x, r)}| dy
\\
&\qquad = \frac{\lambda }{| B(\lambda x, \lambda r)|}\intl_{B(\lambda x, \lambda r)}  | \nabla u(y) - (\nabla u)_{ B(\lambda x, \lambda r)}| dy.
\end{align*}
From this identity we easily deduce   
\begin{equation}
| \nabla u_\lambda |_{ BMO} = \lambda | \nabla  u|_{ BMO}.
\label{A.19}
\end{equation}
We now assume that $ | \nabla u|_{ BMO} >0 $. Otherwise,  since $ \nabla u$ is harmonic and $ u\in L^2(\R^{n})$  would get   
$ u \equiv 0$. Thus,  choosing  
\[
\lambda = \| u\|_{2}^{ \frac{2}{n+2}} | \nabla u|_{ BMO}^{ -\frac{2}{n+2}}, 
\]
by the aid of \eqref{A.18} and \eqref{A.19} we find 
\begin{align*}
\| u\|_{ \infty} &= \| u_{ \lambda }\|_{ \infty} \le c( \| u_\lambda \|_{ 2} + | \nabla u_\lambda |_{ BMO})
=c( \lambda ^{ -\frac{n}{2}} \| u\|_{ 2} + \lambda | \nabla u|_{ BMO})
\\
& =2 c  \| u\|_{2}^{ \frac{2}{n+2}} | \nabla u|_{ BMO}^{ \frac{n}{n+2}}. 
\end{align*}
Whence, \eqref{A.17}.  \hfill \Beweisende  

\vspace{0.5cm}  
In order to estimate the right-hand side of \eqref{A.17} we use the  Calder\'{o}n-Zygmund inequality as follows

\begin{lem}
\label{lemA.10}
For every $ u \in L^2(\R^{2})\cap W^{1,\, 1}_{ loc} (\R^{2})$ with $ \nabla u \in BMO$ the following estimate holds true 
\begin{equation}
| \nabla u|_{ BMO} \le c\Big( | \curl u|_{ BMO} + | \nabla \cdot u|_{ BMO}\Big).   
\label{A.20}
\end{equation}
\end{lem} 

{\bf Proof}:  Denoting the Helmholtz projection  by $ \PP: L^2(\Bbb R^2) \rightarrow L^2(\Bbb R^2)$,  we may write  
$ u = \PP u + (u- \PP u)$.  Clearly, there exists potentials $ \Phi, \Psi \in C^1(\R^{2}) $ with $ \nabla \Phi, \nabla \Psi\in L^2(\R^{2}) $, 
such that    
\[
u = \nabla \Phi + \nabla ^{ \bot} \Psi. 
\]
\vspace{0.5cm}  
Having $ \nabla \cdot u = -\Delta \Psi $, by the Calder\'on-Zygmund inequality, we deduce that  
\[
| \nabla ^2 \Psi |_{ BMO} \le c  | \nabla \cdot u |_{BMO}.
\]
Similarly, observing that $ - \Delta \Psi = \curl u$, once more applying Calder\'on-Zygmund's inequality, we find
\[
| \nabla ^2 \Phi |_{ BMO}  \le c| \curl u |_{ BMO}.  
\] 
From the last two estimates 
we obtain \eqref{A.20}.  \hfill \Beweisende 

\vspace{0.5cm}  
As an immediate consequence of Lemma\,\ref{lemA.9} and Lemma\,\ref{lemA.10}, we get 

\begin{cor}
\label{corA.11}
For every $ u \in L^2(\R^{2})\cap W^{1,\, 1}_{ loc} (\R^{2})$ with $ \nabla u \in BMO$  and 
$ | \nabla \cdot u | \le c | u|$ in $ \R^{2}$ it holds  
\begin{equation}
\|  u\|_{ \infty} \le c \Big(1+\| u\|_{ 2}^{ \frac{1}{2}} | \curl u|_{ BMO}^{ \frac{1}{2}}\Big).  
\label{A.22}
\end{equation}
\end{cor}

\section{Gronwall type iteration  lemma}
\label{sec:-A}
\setcounter{secnum}{\value{section} \setcounter{equation}{0}
\renewcommand{\theequation}{\mbox{B.\arabic{equation}}}}

 \begin{lem}[Iteration lemma]
\label{lemA.4}
Let $ \beta  _m: [t_0, t_1] \rightarrow \R$, $ m\in \N\cup \{ 0\}$ be a sequences of continuous functions. Furthermore let 
$ \alpha, g\in L^1(t_0, t_1)$ with $\alpha \geq 0$. We assume that the following recursive  of integral 
inequality holds true for a constant $ C>0$
\begin{align}
\beta _{ m}(t) \le  Cm + g(t)+   \intl_{t_0}^{t} \alpha(s) \beta _{ m+1}(s)  ds, \quad  m\in \N\cup \{ 0\}.
\label{A.5a}
\end{align}
Furthermore, suppose that there exists a constant $ K>0$ such that 
\begin{equation}
\max_{ t\in [t_0, t_1]}| \beta _m(t)| \le K^m  \quad \forall\,m\in \N.
\label{A.5aa}
\end{equation}
Then the following inequality holds true for all $ t\in [t_0,t_1]$
\begin{equation}
\beta _0 (t)\le   g(t) +C  \intl_{t_0}^{t} \alpha(s)  ds  e^{\intl_{t_0}^{t} \alpha(\tau )}   d\tau 
+ \intl_{t_0}^{t} \alpha(s)  
g(s) e^{\intl_{s}^{t} \alpha(\tau )   d\tau } ds.
\label{A.5b}
\end{equation}
\end{lem}

{\bf Proof}: Iterating \eqref{A.5a} $ m$-times, we see that 
\begin{align}
\beta _0(t) &\le  g(t)+   \intl_{t_0}^{t} \alpha(s_1) \beta _{1}(s_1)  ds_1
\cr
& \le  g(t) +  \intl_{t_0}^{t} \alpha(s_1) (C + g(s_1))   ds_1 + \intl_{t_0}^{t} \intl^{s_1}_{t_0} \alpha(s_1) \alpha(s_2)(2C + g(s_2))   ds_2 ds_1 
\cr
&\qquad + \ldots + \intl_{t_0}^{t} \intl^{s_1}_{t_0} \ldots \intl^{s_{ m-1}}_{t_0} \alpha(s_1) \alpha(s_2) \ldots \alpha(s_m)(Cm + g(s_m))  ds_{ m} \ldots ds_2 ds_1 
\cr
&\qquad + \intl_{t_0}^{t} \intl^{s_1}_{t_0} \ldots \intl^{s_{ m}}_{t_0} \alpha(s_1) \alpha(s_2) \ldots \alpha(s_{ m+1}) \beta _{ m+1}(s_{ m+1})  
ds_{ m+1} \ldots ds_2 ds_1 
\cr
&= g(t) +    \intl_{t_0}^{t} \alpha(s) (C+ g(s))  ds+   \sum_{k=2}^{m} I_k + J_m,
\label{A.14}
\end{align} 
 where 
\[
I_k = \intl_{t_0}^{t} \intl^{s_1}_{t_0} \ldots \intl^{s_{ k-1}}_{t_0} \alpha(s_1) \alpha(s_2) \ldots \alpha(s_k)(Ck + g(s_k))  ds_{ k} \ldots ds_2 ds_1, \quad  k=2, \ldots, m.  
\]
With the help of Fubini's theorem   we compute 
\begin{align*}
 I_k = \intl_{t_0}^{t} \intl^{t}_{s_k} \ldots \intl^{s_{ k-2}}_{s_k} \alpha(s_1) \alpha(s_2) \ldots \alpha(s_{ k-1}) ds_{ k-1}  \ldots ds_2 ds_1 \alpha(s_{ k})(Ck + g(s_k) )ds_{ k}
\end{align*}
Using the following identity
\begin{equation}
 \intl_{s}^{t}\intl^{s_1}_{s} \ldots \intl^{s_{ k-2}}_{s}   
 \prod_{ j=1}^{ k-1} \alpha (s_j)   d s_{ k-1} \ldots d s_1
 = \frac{1}{(k-1)!}  \bigg(\intl_{s}^{t} \alpha (\tau ) d\tau \bigg)^{ k-1}, 
\label{A.2}
\end{equation}
we obtain 
 \begin{align*}
I_k
&=  \intl_{t_0}^{t}   \alpha(s)(Ck + g(s)) \frac{1}{(k-1)!}  \bigg(\intl_{s}^{t} \alpha(\tau )   d\tau \bigg)^{ k-1} ds
\\
&= -\frac{C}{(k-1)!} \intl_{t_0} ^t  \frac{d}{ds}\bigg(\intl_{s}^{t} \alpha(\tau )   d\tau \bigg)^{ k} ds +
\intl_{t_0}^{t}   \alpha(s)g(s)\frac{1}{(k-1)!}  \bigg(\intl_{s}^{t} \alpha(\tau )   d\tau \bigg)^{ k-1} ds\\
&=\frac{C}{(k-1)!} \bigg(\intl_{t_0}^{t} \alpha(\tau )   d\tau \bigg)^{ k}  +
\intl_{t_0}^{t}   \alpha(s)g(s)\frac{1}{(k-1)!}  \bigg(\intl_{s}^{t} \alpha(\tau )   d\tau \bigg)^{ k-1} ds.
\end{align*}
Furthermore, by our assumption on $ \{ \beta _m\}$ along with \eqref{A.2} with $ k= m+1$  we see that 
\begin{align*}
 | J_m| &= \bigg|\intl_{t_0}^{t}   \beta _{ m+1}(s) \frac{1}{m!}  \bigg(\intl_{s}^{t} \alpha(\tau )   d\tau \bigg)^{m} ds\bigg| 
\\
& \le K (t_1-t_0)  \frac{1}{m!}\bigg(K \intl_{t_0}^{t_1}  \alpha(\tau )   d\tau \bigg)^m \rightarrow 0 \quad  \text{ as}\quad 
m \rightarrow +\infty.
\end{align*}
Therefore, letting $ m \rightarrow \infty$ in the right-hand side of \eqref{A.14}, we arrive at  
 \begin{align*}
\beta _0 (t)\le   g(t) +C  \intl_{t_0}^{t} \alpha(s)  ds  e^{\intl_{t_0}^{t} \alpha(\tau )   d\tau }
+ \intl_{t_0}^{t} \alpha(s)  
g(s) e^{\intl_{s}^{t} \alpha(\tau )   d\tau } ds.
\end{align*}
Whence, \eqref{A.5b} follows.  
\hfill \Beweisende  

\section{Continuity for solutions to the transport equation}
\label{sec:-C}
\setcounter{secnum}{\value{section} \setcounter{equation}{0}
\renewcommand{\theequation}{\mbox{C.\arabic{equation}}}}

\begin{lem}
\label{lemC.1}
Given  $ f\in L^q(\R^{n}\times (t_0, 0))$, $ 1< q< +\infty$,  and $ u\in L^\infty(t_0,0; W^{1,\, \infty}(\R^{n}) )
\cap C(\R^{n}\times [t_0, 0])$, let  $ h\in  L^q(\R^{n}\times (t_0,0))$ be a weak solution to  
\begin{equation}
\partial _t h + u\cdot \nabla h = f \quad  \text{ in}\quad  \R^{n}\times (t_0, 0). 
\label{C.1}
\end{equation}
Then,  $ h\in C([t_0, 0 ]; L^q(\R^{n}))$, and it holds for all $ t\in [t_0, 0]$
 \begin{align}
\| h(t)\|_{ q}^q &= \| h(t_0)\|_{ q}^q  + \intl_{t_0}^{t}\intl_{\Bbb R^n} \nabla \cdot  u | h|^{ q}  dxds   + q \intl_{t_0}^{t}\intl_{\Bbb R^n} f h | h|^{ q-2}  dxds.
\label{C.2}
\end{align}

\end{lem}

{\bf Proof}:  
By $ \eta \in C^\infty_c(B(1))$ we denote the usual mollifying kernel. We set 
$ \eta _\var(x) = \var ^{ -n}\eta (\var ^{ -1} x)$, $ 0< \var < +\infty$,  and define for $ v\in L^1_{ loc}(\R^{n})$ 
\[
v_\var (x) = \intl_{ \R^{n}}v(x-y) \eta _\var (y) dy= \intl_{ \R^{n}}v(y) \eta _\var (x-y) dy, \quad  x\in \R^{n}.  
\]
Testing \eqref{C.1} with $ \eta_\var  (x - \cdot )$,   we see that  $ h_\var $ solves 
\begin{equation}
\partial _t h_\var + \nabla \cdot (uh)_{ \var } = f_\var  + (\nabla \cdot u h)_\var.   
\label{C.4}
\end{equation}
Clearly, $ h_\var \in L^q(t_0,0; L^q(\R^{n}))$ with $ \partial _t h_\var \in L^q(t_0,0; L^q(\R^{n}))$. Eventually, redefining 
$ h_\var(t) $ on a set of measure zero in $ [t_0,0]$, we may assume that $ h_\var \in C([t_0, 0]; L^q(\R^{n}))$.  

We estimate for $ (x,t)\in \R^{n}\times [t_0, 0]$
\begin{align*}
&\nabla \cdot (uh)_{ \var }(x, t) - (u \cdot \nabla h_{ \var })(x, t)  
\\
&= \intl_{ \R^{n}} u(x-y, t) h(x-y,t)\cdot \nabla \eta _\var (y)  - u(x,t)h(x-y,t)\cdot \nabla \eta _\var(y)   dy
\\
&= \intl_{ \R^{n}} (u(x-y,t) - u(x,t)) (h(x-y,t)- h(x,t))\cdot \nabla \eta _\var (y)    dy.
\end{align*}
Accordingly,
\begin{equation}
| \nabla \cdot (uh)_{ \var }(x,t) -   \nabla \cdot (uh_{ \var })(x, t)  | \le c \| \nabla u(t)\|_{ \infty} \intmw_{B(\var )} | h(x-y,t)-h(x,t)| dy 
\label{C.5}
\end{equation}
Multiplying this inequality by $ | h_\var (x,t)|^{ q-1}$, integrating over $ \R^{n}$ with respect to $ x$,  and applying H\"older's inequality, we find 
\begin{align}
&\intl_{\R^{n}}| \nabla \cdot (uh)_{ \var }(x,t) -(u\cdot \nabla h_{ \var })(x, t)  | | h_\var(x,t) |^{ q-1}dx
\cr
&\quad \le c \| \nabla u(t)\|_{ \infty} \intl_{\R^{n}}\intmw_{B(\var )} | h(x-y,t)-h(x,t)| | h_\var (x,t)|^{ q-1} dydx 
\cr
&\quad \le c \| \nabla u(t)\|_{ \infty} \intmw_{B(\var )}\intl_{\R^{n}} | h(x-y,t)-h(x,t)| | h_\var (x,t)|^{ q-1} dx dy
\cr
&\quad \le c \| \nabla u(t)\|_{ \infty} \intmw_{B(\var )} \| h(\cdot -y,t) -h(t)\|_{ q}dy \| h\|_{ q}^{ q-1} 
\le c \| \nabla u(t)\|_{ \infty} \| h(t)\|_q^q. 
\label{C.6a}
\end{align}
Multiplying \eqref{C.4}  by $ qh_\var | h_\var |^{ q-2}$, and applying integration by parts, together with \eqref{C.6a} we see 
that for almost every $ t\in (t_0, 0)$
\begin{align*}
\| h_\var(t) \|_q^q &= \| h_\var (t_0)\|_{ q}^q-q \intl_{t_0}^{t} \intl_{\R^{n}}
\Big(\nabla \cdot (uh)_{ \var } -  u\cdot \nabla h_{ \var }\Big)  h_\var| h_\var|^{ q-2}  dx ds
\\
&\qquad  +q \intl_{t_0}^{t} \intl_{\R^{n}}
(\nabla \cdot u h)_\var   h_\var| h_\var |^{ q-2} dx ds + \intl_{t_0}^{t} \intl_{\R^{n}}
\nabla \cdot u | h_\var |^q  dx ds 
\\
& \qquad + \intl_{t_0}^{t} \intl_{\R^{n}} f_\var h_\var| h_\var|^{ q-2} dxds
\\
& \le  \| h (t_0)\|_{ q}^q + c (1+\| \nabla u\|_{ \infty})\| h\|_{ q}^q
 + \| f\|_{ q}^q. 
\end{align*}
This yields, $ h\in L^\infty(t_0, 0 ; L^q(\R^{n}))$. On the other hand, from \eqref{C.1} we deduce that $ \partial _t h
\in L^q(t_0, 0; W^{-1,\, q}(\R^{n}))$. Thus, eventually redefining $ h(t)$ on a set of measure zero, we get $ h\in C_w([t_0, 0]; L^q(\R^{n}))$.    

Let $ t\in (t_0,0]$. Multiplying \eqref{C.4} by $ qh_\var | h_\var |^{ q-2}$  and integrating over $ \R^{n}\times (t_0, t)$,  we find 
\begin{align}
\| h_\var(t) \|_q^q &= \| h_\var (t_0)\|_{ q}^q - q \intl_{t_0}^{t} \intl_{\R^{n}}
\Big(\nabla \cdot (uh)_{ \var } - (u \cdot \nabla h_{ \var })\Big)  h_\var | h_\var |^{ q-2}  dx ds
\cr
&\qquad +q \intl_{t_0}^{t} \intl_{\R^{n}}
(\nabla \cdot u h)_\var h_\var  | h_\var|^{ q-2} dx ds + \intl_{t_0}^{t} \intl_{\R^{n}}
\nabla \cdot u | h_\var |^q  dx ds
\cr
& \qquad + q \intl_{t_0}^{t} \intl_{\R^{n}} f_\var h_\var| h_\var |^{ q-2}dxds.
\label{C.5a}
\end{align}

We claim that 
\begin{align}
\begin{cases}
{\D -\intl_{t_0}^{t} \intl_{\R^{n}}}
\Big(\nabla \cdot (uh)_{ \var }- (u\cdot \nabla h_{ \var })\Big)  h_\var | h_\var |^{ q-2}  dx ds 
\\
\qquad  +  {\D\intl_{t_0}^{t} \intl_{\R^{n}}}
(\nabla \cdot u h)_\var  h_\var  | h_\var |^{ q-2} dx ds
\rightarrow 0
\quad \text{{\it as}}\quad  \var \rightarrow 0. 
\end{cases}
\label{C.5b}
\end{align}

To see this, we first verify that 
 \begin{equation}
-\nabla \cdot (uh)_{ \var } + u\cdot \nabla h_{ \var }+(\nabla \cdot u h)_\var :=\sigma_\var  \rightarrow 0 \quad \text{{\it  weakly in }}\quad  L^q(\R^{n}\times (t_0, 0))
\label{C.6}
\end{equation}
 as $\var \rightarrow 0$.
Indeed, thanks to \eqref{C.5}  we see that $\sigma_\var \in L^q(\R^{n}\times (t_0, 0))$. 
Let $ \{ \var _k\}$ be a sequence of positive numbers with $\var_k\to  0$ as $ k \rightarrow +\infty$. Eventually passing to a subsequence, we get $ \sigma \in L^q(\R^{n}\times (t_0, 0))$ such that 
\[
\sigma_{\var_k}
\rightarrow \sigma \quad \text{ {\it weakly in} }\quad  L^q(\R^{n}\times (t_0, 0))\quad
\text{{\it as}}\quad k \rightarrow +\infty.
\]
On the other hand, we get for all $ \varphi \in C^\infty_c(\R^{n}\times (t_0, 0))$ 
 \begin{align*}
&\intl_{t_0}^{0} \intl_{\R^{n}}
\Big(-\nabla \cdot (uh)_{ \var_k } + u\cdot \nabla h_{ \var_k }+(\nabla \cdot u h)_{ \var _k}  \Big)  \varphi dx ds 
\\
&\,   = \intl_{t_0}^{0} \intl_{\R^{n}}\Big( (uh)_{ \var_k } -  (uh_{ \var_k })\Big)\cdot   \nabla \varphi  dx ds 
  + \intl_{t_0}^{0} \intl_{\R^{n}}\Big( -\nabla \cdot u h_{ \var_k } + (\nabla \cdot u h)_{ \var_k }\Big)\varphi  dx ds 
\\
&\quad  \rightarrow 0\quad \text{{\it as}}\quad  k \rightarrow \infty. 
\end{align*}
This provides us with 
\[
\intl_{t_0}^{0} \intl_{\R^{n}} \sigma  \varphi dxds=0\quad  \forall\,\varphi \in C^\infty_c(\R^{n}\times (t_0, 0)). 
\]
 Accordingly, $ \sigma =0$.  Whence, \eqref{C.6}.  
 
Next,  by means of Lebesgue's theorem of dominated convergence, we see that 
\begin{equation}
h_{ \var } | h_{ \var} |^{ q-2}-h | h |^{ q-2}  \rightarrow 0\quad  \text{ in}\quad  L^{ q'}(\R^{n}\times (t_0, 0))\quad  \text{ as}\quad  
\var  \rightarrow 0. 
\label{C.7}
\end{equation}
Clearly, 
\begin{align*}
&-\intl_{t_0}^{t} \intl_{\R^{n}}
\Big(\nabla \cdot (uh)_{ \var } - \nabla \cdot (uh_{ \var })\Big)  h_{ \var }| h_{ \var }|^{ q-2}  dx ds 
+\intl_{t_0}^{t} \intl_{\R^{n}}
(\nabla \cdot u h)_\var  h_\var  | h_\var |^{ q-2} dx ds
\\
&\quad = \intl_{t_0}^{0} \intl_{\R^{n}}
\Big(-\nabla \cdot (uh)_{ \var }+ u\cdot \nabla h_{ \var }+ (\nabla \cdot u h)_\var \Big)  \Big(h_{ \var } | h_{ \var } |^{ q-2} - h| h|^{ q-2}\Big) dx ds 
\\
&\qquad +\intl_{t_0}^{0} \intl_{\R^{n}}
\Big(-\nabla \cdot (uh)_{ \var }+ u\cdot \nabla h_{ \var }+ (\nabla \cdot u h)_\var \Big)  h | h |^{ q-2}  dx ds. 
\end{align*} 

The first integral on the right-hand side of the above identity  converges to $ 0$ in view of \eqref{C.7}, while the second integral converges to $ 0$ according 
to \eqref{C.6}. Whence, \eqref{C.5b}.     
 
Thus, thanks to  \eqref{C.7}, \eqref{C.5a} and the fact that $ h_{ \var } (t) \rightarrow h$ in $ L^q(\R^{n})$ for all $ t\in [t_0, 0]$ 
and $ f_\var \rightarrow f $ in $ L^q(\R^{n}\times (t_0,0))$ as $ \var \rightarrow 0$, we   
are in a position to pass $ \var \rightarrow 0$ in \eqref{C.5a}. This yields the following identity 
\begin{align}
\| h(t) \|_q^q = \| h (t_0)\|_{ q}^q + \intl_{t_0}^{t} \intl_{\R^{n}} \nabla \cdot u | h |^{ q} dxds + q \intl_{t_0}^{t} \intl_{\R^{n}} f h | h|^{ q-2}dxds.
\label{C.8}
\end{align}
In particular, $ \| h(t) \|_q^q \in C([t_0,0])$, which together with $ h\in C_w([t_0,0; L^q(\R^{n})])$ yields $ h\in C([t_0,0]; L^q(\R^{n}))$. 
This completes the proof of the lemma.  \hfill \Beweisende

$$\mbox{\bf Acknowledgements}$$
Chae was partially supported by NRF grants 2016R1A2B3011647, while Wolf has been supported 
supported by NRF grants 2017R1E1A1A01074536.  
  The authors declare that they have no conflict of interest.
  \bibliographystyle{siam} 

\end{document}